
\def\autori{G.\ DAL MASO and I.V.\ SKRYPNIK}
\def\titolo{Capacity theory for monotone operators}




\font\sixrm=cmr6
\newcount\tagno \tagno=0		        
\newcount\thmno	\thmno=0	         	
\newcount\bibno	\bibno=0			
\newcount\chapno\chapno=0                       
\newcount\verno            
\newif\ifproofmode
\proofmodetrue
\newif\ifwanted
\wantedfalse
\newif\ifindexed
\indexedfalse

\def\ifundefined#1{\expandafter\ifx\csname+#1\endcsname\relax}

\def\Wanted#1{\ifundefined{#1} \wantedtrue
\immediate\write0{Wanted
#1
\the\chapno.\the\thmno}\fi}

\def\Increase#1{{\global\advance#1 by 1}}

\def\Assign#1#2{\immediate
\write1{\noexpand\expandafter\noexpand\def
 \noexpand\csname+#1\endcsname{#2}}\relax
 \global\expandafter\edef\csname+#1\endcsname{#2}}

\def\pAssign#1#2{\write1{\noexpand\expandafter\noexpand\def
 \noexpand\csname+#1\endcsname{#2}}}

\def\lPut#1{\ifproofmode\llap{\hbox{\sixrm #1\ \ \ }}\fi}
\def\rPut#1{\ifproofmode$^{\hbox{\sixrm #1}}$\fi}



\def\chp#1{\global\tagno=0\global\thmno=0\Increase\chapno
\Assign{#1}
{\the\chapno}{\lPut{#1}\the\chapno}}


\def\thm#1{\Increase\thmno
\Assign{#1}
{\the\chapno.\the\thmno}\the\chapno.\the\thmno\rPut{#1}}


\def\frm#1{\Increase\tagno
\Assign{#1}{\the\chapno.\the\tagno}\lPut{#1}
{\the\chapno.\the\tagno}}


\def\bib#1{\Increase\bibno
\Assign{#1}{\the\bibno}\lPut{#1}{\the\bibno}}


\def\pgp#1{\pAssign{#1/}{\the\pageno}}


\def\ix#1#2#3{\pAssign{#2}{\the\pageno}
\immediate\write#1{\noexpand\idxitem{#3}
{\noexpand\csname+#2\endcsname}}}


\def\rf#1{\Wanted{#1}\csname+#1\endcsname\relax\rPut {#1}}


\def\rfp#1{\Wanted{#1}\csname+#1/\endcsname\relax\rPut{#1}}

\input \jobname.refs
\Increase\verno
\immediate\openout1=\jobname.aux

\immediate\write1{\noexpand\verno=\the\verno}

\ifindexed
\immediate\openout2=\jobname.idx
\immediate\openout3=\jobname.sym
\fi


\font\twelverm=cmr12
\font\twelvei=cmmi12
\font\twelvesy=cmsy10
\font\twelvebf=cmbx12
\font\twelvett=cmtt12
\font\twelveit=cmti12
\font\twelvesl=cmsl12

\font\ninerm=cmr9
\font\ninei=cmmi9
\font\ninesy=cmsy9
\font\ninebf=cmbx9
\font\ninett=cmtt9
\font\nineit=cmti9
\font\ninesl=cmsl9

\font\eightrm=cmr8
\font\eighti=cmmi8
\font\eightsy=cmsy8
\font\eightbf=cmbx8
\font\eighttt=cmtt8
\font\eightit=cmti8
\font\eightsl=cmsl8

\font\sixrm=cmr6
\font\sixi=cmmi6
\font\sixsy=cmsy6
\font\sixbf=cmbx6

\catcode`@=11 
\newskip\ttglue

\def\twelvepoint{\def\rm{\fam0\twelverm}
\textfont0=\twelverm  \scriptfont0=\ninerm
\scriptscriptfont0=\sevenrm
\textfont1=\twelvei  \scriptfont1=\ninei  \scriptscriptfont1=\seveni
\textfont2=\twelvesy  \scriptfont2=\ninesy
\scriptscriptfont2=\sevensy
\textfont3=\tenex  \scriptfont3=\tenex  \scriptscriptfont3=\tenex
\textfont\itfam=\twelveit  \def\it{\fam\itfam\twelveit}%
\textfont\slfam=\twelvesl  \def\sl{\fam\slfam\twelvesl}%
\textfont\ttfam=\twelvett  \def\tt{\fam\ttfam\twelvett}%
\textfont\bffam=\twelvebf  \scriptfont\bffam=\ninebf
\scriptscriptfont\bffam=\sevenbf  \def\bf{\fam\bffam\twelvebf}%
\tt  \ttglue=.5em plus.25em minus.15em
\normalbaselineskip=15pt
\setbox\strutbox=\hbox{\vrule height10pt depth5pt width0pt}%
\let\sc=\tenrm  \let\big=\twelvebig  \normalbaselines\rm}

\def\tenpoint{\def\rm{\fam0\tenrm}
\textfont0=\tenrm  \scriptfont0=\sevenrm  \scriptscriptfont0=\fiverm
\textfont1=\teni  \scriptfont1=\seveni  \scriptscriptfont1=\fivei
\textfont2=\tensy  \scriptfont2=\sevensy  \scriptscriptfont2=\fivesy
\textfont3=\tenex  \scriptfont3=\tenex  \scriptscriptfont3=\tenex
\textfont\itfam=\tenit  \def\it{\fam\itfam\tenit}%
\textfont\slfam=\tensl  \def\sl{\fam\slfam\tensl}%
\textfont\ttfam=\tentt  \def\tt{\fam\ttfam\tentt}%
\textfont\bffam=\tenbf  \scriptfont\bffam=\sevenbf
\scriptscriptfont\bffam=\fivebf  \def\bf{\fam\bffam\tenbf}%
\tt  \ttglue=.5em plus.25em minus.15em
\normalbaselineskip=12pt
\setbox\strutbox=\hbox{\vrule height8.5pt depth3.5pt width0pt}%
\let\sc=\eightrm  \let\big=\tenbig  \normalbaselines\rm}

\def\ninepoint{\def\rm{\fam0\ninerm}
\textfont0=\ninerm  \scriptfont0=\sixrm  \scriptscriptfont0=\fiverm
\textfont1=\ninei  \scriptfont1=\sixi  \scriptscriptfont1=\fivei
\textfont2=\ninesy  \scriptfont2=\sixsy  \scriptscriptfont2=\fivesy
\textfont3=\tenex  \scriptfont3=\tenex  \scriptscriptfont3=\tenex
\textfont\itfam=\nineit  \def\it{\fam\itfam\nineit}%
\textfont\slfam=\ninesl  \def\sl{\fam\slfam\ninesl}%
\textfont\ttfam=\ninett  \def\tt{\fam\ttfam\ninett}%
\textfont\bffam=\ninebf  \scriptfont\bffam=\sixbf
\scriptscriptfont\bffam=\fivebf  \def\bf{\fam\bffam\ninebf}%
\tt  \ttglue=.5em plus.25em minus.15em
\normalbaselineskip=11pt
\setbox\strutbox=\hbox{\vrule height8pt depth3pt width0pt}%
\let\sc=\sevenrm  \let\big=\ninebig  \normalbaselines\rm}

\def\eightpoint{\def\rm{\fam0\eightrm}
\textfont0=\eightrm  \scriptfont0=\sixrm  \scriptscriptfont0=\fiverm
\textfont1=\eighti  \scriptfont1=\sixi  \scriptscriptfont1=\fivei
\textfont2=\eightsy  \scriptfont2=\sixsy  \scriptscriptfont2=\fivesy
\textfont3=\tenex  \scriptfont3=\tenex  \scriptscriptfont3=\tenex
\textfont\itfam=\eightit  \def\it{\fam\itfam\eightit}%
\textfont\slfam=\eightsl  \def\sl{\fam\slfam\eightsl}%
\textfont\ttfam=\eighttt  \def\tt{\fam\ttfam\eighttt}%
\textfont\bffam=\eightbf  \scriptfont\bffam=\sixbf
\scriptscriptfont\bffam=\fivebf  \def\bf{\fam\bffam\eightbf}%
\tt  \ttglue=.5em plus.25em minus.15em
\normalbaselineskip=9pt
\setbox\strutbox=\hbox{\vrule height7pt depth2pt width0pt}%
\let\sc=\sixrm  \let\big=\eightbig  \normalbaselines\rm}

\def\twelvebig#1{{\hbox{$\textfont0=\twelverm\textfont2=\twelvesy
	\left#1\vbox to10pt{}\right.\n@space$}}}
\def\tenbig#1{{\hbox{$\left#1\vbox to8.5pt{}\right.\n@space$}}}
\def\ninebig#1{{\hbox{$\textfont0=\tenrm\textfont2=\tensy
	\left#1\vbox to7.25pt{}\right.\n@space$}}}
\def\eightbig#1{{\hbox{$\textfont0=\ninerm\textfont2=\ninesy
	\left#1\vbox to6.5pt{}\right.\n@space$}}}

\def\displayliness#1{\null\,\vcenter{\openup1\jot \m@th
  \ialign{\strut\hfil$\displaystyle{##}$\hfil
    \crcr#1\crcr}}\,}
\def\displaylinesno#1{\displ@y \tabskip=\centering
   \halign to\displaywidth{ \hfil$\@lign \displaystyle{##}$ \hfil
    	\tabskip=\centering
     &\llap{$\@lign##$}\tabskip=0pt \crcr#1\crcr}}
\def\ldisplaylinesno#1{\displ@y \tabskip=\centering
   \halign to\displaywidth{ \hfil$\@lign \displaystyle{##}$\hfil
    	\tabskip=\centering
     &\kern-\displaywidth
     \rlap{$\@lign##$}\tabskip=\displaywidth \crcr#1\crcr}}

\catcode`@=12 

\def\parag#1#2{\goodbreak\bigskip\bigskip\noindent
                   {\bf #1.\ \ #2}
                   \nobreak\bigskip}
\def\intro#1{\goodbreak\bigskip\bigskip\goodbreak\noindent
                   {\bf #1}\nobreak\bigskip\nobreak}
\long\def\th#1#2{\goodbreak\bigskip\noindent
                {\bf Theorem #1.\ \ \it #2}}
\long\def\lemma#1#2{\goodbreak\bigskip\noindent
                {\bf Lemma #1.\ \ \it #2}}
\long\def\prop#1#2{\goodbreak\bigskip\noindent
                  {\bf Proposition #1.\ \ \it #2}}
\long\def\cor#1#2{\goodbreak\bigskip\noindent
                {\bf Corollary #1.\ \ \it #2}}
\long\def\defin#1#2{\goodbreak\bigskip\noindent
                  {\bf Definition #1.\ \ \rm #2}}
\long\def\rem#1#2{\goodbreak\bigskip\noindent
                 {\bf Remark #1.\ \ \rm #2}}

\def\spazio{\bigskip}
\def\negbigskip{\vskip-\bigskipamount}

\def\proof{\vskip.4cm\noindent{\it Proof.\ \ }}

\mathchardef\emptyset="001F
\mathchardef\hyphen="002D


\def\meas{{\rm meas}}
\def\cp{C_p}
\def\ca{C_A}
\def\cas{{\hat C}_A}
\def\supp{{\rm supp}}
\def\div{{\rm div}}

\def\dirpb#1#2#3#4{\cases{#1\,,&\cr\cr #2
\qquad \hbox{in\ }\,#3\,#4&\cr}}


\def\Eq{\,=\,}

\def\Ge{\,\ge\,}
\def\Le{\,\le\,}

\def \r{{\bf R}}
\def\rn{\r^n}

\def \n{{\bf N}}

\def\interior{\mathaccent'27}

\def\pinfty{{+}\infty}
\def\finite{<\pinfty}

\def\e{\varepsilon}

\def\Om{\Omega}



\def\cinftyo#1{C^\infty_0(#1)}
\def\lone#1{L^1_{}(#1)}
\def\lp#1{L^p_{}(#1)}
\def\lpn#1{L^p_{}(#1,\rn)}
\def\lq#1{L^q_{}(#1)}
\def\linfty#1{L^\infty(#1)}
\def\qe{$\cp$-q.e.\ in}
\def\ae{a.e.\ in}
\def\wp#1{W^{1,p}_{}(#1)}
\def\wo#1{W^{1,p}_0(#1)}
\def\wm#1{W^{-1,q}_{}(#1)}

\def\into{\int_\Omega}

\def\A#1{A#1}
\def\a#1{a(x,#1)}
\def\aa#1#2{\bigl(\a{#1},#2\bigr)}
\def\aaa#1#2#3{\bigl(\a{#1}-\a{#2},#3\bigr)}
\def\iaa#1#2#3{\int_{#1} \aa{#2}{#3}\,dx}
\def\iaaa#1#2#3#4{\int_{#1} \aaa{#2}{#3}{#4}\,dx}

\def\sqr#1#2{\vbox{
   \hrule height .#2pt
   \hbox{\vrule width .#2pt height #1pt \kern #1pt
      \vrule width .#2pt}
   \hrule height .#2pt }}
\def\square{\sqr74}

\def\endproof{{\unskip\nobreak\hfill
\hskip2em\hbox{}\nobreak\hfill $\square$ \goodbreak
\parfillskip=0pt  \finalhyphendemerits=0}}

\def\rightheadline{\eightpoint\hfil\titolo
\hfil\tenrm\folio}
\def\leftheadline{\tenrm\folio\hfil\eightpoint
\autori \hfil}
\def\zeroheadline{\hfill}
%
\headline={\ifnum\pageno=0 \zeroheadline
\else\ifodd\pageno\rightheadline
\else\leftheadline\fi\fi}

\nopagenumbers
\magnification=1200
\baselineskip=15pt
\hfuzz=2pt
\parindent=2em
\mathsurround=1pt
\tolerance=1000

\pageno=0
\hsize 14truecm
\vsize 25truecm
\hoffset=0.8truecm
\voffset=-1.55truecm

\null
\vskip 2.8truecm
{\twelvepoint
\baselineskip=1.7\baselineskip
\centerline{\bf CAPACITY THEORY FOR}
\centerline{\bf MONOTONE OPERATORS}
}
\vskip2truecm
\centerline{Gianni DAL MASO ($^1$)}
\medskip
\centerline{Igor V. SKRYPNIK ($^2$)}
\vfil

{\eightpoint
\baselineskip=1.2\baselineskip
\centerline{\bf Abstract}
\bigskip
\noindent If $\A{u}=-\div\bigl(\a{Du}\bigr)$ is a monotone operator
defined on the Sobolev space $\wp{\rn}$, ${1<p\finite}$, with $\a{0}=0$
for a.e.\ $x\in\rn$, the capacity $\ca(E,F)$ relative to~$A$ can be
defined for every pair $(E,F)$ of bounded sets in~$\rn$ with
$E\subset F$. We prove that $\ca(E,F)$ is increasing and countably
subadditive with respect to~$E$ and decreasing with respect to~$F$.
Moreover we investigate the continuity properties of $\ca(E,F)$ with
respect to~$E$ and~$F$.
\par
\vfil
\item{$(^1)$} SISSA, Via Beirut 4, 34013 Trieste, Italy
\item{} e-mail: dalmaso@tsmi19.sissa.it
\smallskip
\item{$(^2)$} Institute of Applied Mathematics and Mechanics,
\item{} Academy of Sciences of Ukraine, R. Luxemburg St. 74,
\item{} 340114 Donetsk, Ukraine
\item{} e-mail: skrypnik@iamm.ac.msk.su

\par
}
\vskip 1truecm
\centerline {Ref. S.I.S.S.A.
6/95/M (January 95)}
\vskip 1truecm
\eject
%
%
%
\topskip=25pt 
\vsize 22.5truecm
\hsize 16.2truecm
\hoffset=0truecm
\voffset=0.5truecm

\proofmodefalse 

\intro{Introduction}
Let $A\colon\wp{\rn}\to\wm{\rn}$, ${1<p\finite}$, ${1/p+1/q=1}$, be a
monotone operator of the form
$$
\A{u}=-\div\bigl(\a{Du}\bigr)\,,
\leqno(\frm{monop})
$$
where $a\colon\rn{\times}\rn\to\rn$ is a Carath\'eodory function
which satisfies the usual monotonicity, coerciveness, and growth
conditions  (see~(\rf{mono}), (\rf{coerc}), (\rf{bound}) below), and
$\a{0}=0$ for a.e.\ $x\in\rn$.

If~$E$ and~$F$ are bounded sets in~$\rn$, with~$E$ closed, $F$
open, and $E\subset F$, the capacity of~$E$ in~$F$ relative to the
operator~$A$ is defined  as
$$
\ca(E,F)\Eq \iaa{F\setminus E}{Du}{Du}\,,
\leqno(\frm{cap0})
$$
where~$u$, the $\ca$-potential of~$E$ in~$F$, is the weak solution of the
Dirichlet problem
$$
\A{u}=0\,\ \hbox{in\ }\,F\setminus E\,, \qquad
u=1\,\ \hbox{in\ }\,\partial E\,,
\,\ u=0\,\ \hbox{in\ }\,\partial F\,.
\leqno(\frm{dir0})
$$
This definition is extended to arbitrary bounded
sets by giving a  suitable meaning to problem~(\rf{dir0}) when
${F\setminus E}$ is not  open (Definition~\rf{cap}).

The purpose of this paper is to prove the main properties of the set
function~$\ca$. In particular we prove that $\ca(E,F)$ is increasing with
respect to~$E$ (Theorem~\rf{monoE}) and decreasing with respect
to~$F$ (Theorem~\rf{monoF}). Moreover, we show that $\ca(\cdot,F)$ is
continuous along all increasing sequences of sets
(Theorem~\rf{increasingE}) and along all decreasing sequences of
closed sets contained in the interior of~$F$
(Theorem~\rf{decreasingE}), while $\ca(E,\cdot)$ is continuous along all
decreasing sequences of sets  (Theorem~\rf{decreasingF}) and along all
increasing sequences of  open sets containing the closure of~$E$
(Theorem~\rf{increasingF}). These results allow us to show that
$$
\ca(E,F)\Eq
\sup\{\ca(K,U): K\,\ \hbox{compact,\ } \, K\subset E\,,
\,\ U\,\ \hbox{bounded\ and\ open,\ } \, U\supset F \}
$$
when~$E$ and~$F$ are bounded Borel sets (Theorem~\rf{choquet}), and
to prove that  $\ca(\cdot,F)$ is countably subadditive
(Theorem~\rf{countsub}).

Finally, we introduce the capacity $\ca(E,F,s)$ with respect to a constant
${s\in\r}$ by replacing the condition ${u=1}$ in $\partial E$
which appears in~(\rf{dir0}) with the condition ${u=s}$ in $\partial E$
(Definition~\rf{caps}). We prove that the function ${1\over s}\ca(E,F,s)$
is continuous and increasing with respect to~$s$ (Theorems~\rf{caps15}
and~\rf{monos}).

When $\A{u}=-\div\bigl(|Du|^{p-2}Du\bigr)$, the capacity~$\ca$
coincides with  the usual capacity $\cp$ associated with the Sobolev
space $\wp{\rn}$  (see Section~\rf{notation}), for which the above
mentioned properties are well known and can be obtained easily by
using the fact that~(\rf{dir0}) is the Euler  equation of a suitable
minimum problem, and thus $\cp(E,F)$ can be defined  equivalently as
the infimum of $\int_F |Dv|^pdx$ over the set of all functions~$v$ in
$\wo{F}$ such that $v\ge 1$ in a neighbourhood of~$E$. For a monotone
operator of the  form~(\rf{monop}) problem~(\rf{dir0}) is, in  general,
not equivalent to a minimum  problem, and this fact forces us
to develop a completely new proof.

When the operator~$A$ is linear, the capacity~$\ca$ was introduced
in~[\rf{STA}], but, to our knowledge, the properties considered above
have been
established only in~[\rf{DM-GAR-cap}]. The proof avoids auxiliary
minimum
problems, but involves the adjoint operator~$A^*$ in an essential way,
and, therefore, it can not be adapted to the monotone case.

Our proof is based on an estimate of the $\ca$-potentials
(Lemmas~\rf{comp1} and~\rf{comp2}), which follows from a standard
comparison argument (Theorem~\rf{comparison}). The main
new idea is to deduce the inequalities for the capacity~$\ca$ from the
corresponding inequalities for the $\ca$-potentials. The tools used in
this approach are the notion of $\ca$-distribution
(Theorem~\rf{capdistr1}, Definition~\rf{capdistr}, and
Proposition~\rf{cd}) and a technical
lemma  which allows us, under very special conditions, to deduce the
inequality  $\A{u_1}\ge\A{u_2}$ from the inequality ${u_1\le u_2}$
(Lemma~\rf{fundamental}).

For a complete treatment of the problem, we consider also the case
when the operator~$A$ is not strictly monotone, and thus~(\rf{dir0})
may have more than one solution. We prove that in this case the
capacity
$\ca(E,F)$ defined by~(\rf{cap0}) does not depend on the choice of the
$\ca$-potential~$u$. The proof is based on a careful investigation,
developed in Section~\rf{estimates}, of the properties of the set of all
solutions of problem~(\rf{dir0}).

Under some natural assumptions on~$a$ the capacities $\ca$ and
$\cp$ are equivalent (Remark~\rf{cap2}), i.e., there exist two constants
${\alpha>0}$ and ${\beta>0}$ such that
$$
\alpha\, \cp(E,F) \le \ca(E,F) \le \beta\, \cp(E,F) \,.
$$
Therefore the precise behaviour of the capacity~$\ca$ is not
important in all those problems where it is enough to obtain just an
estimate of $\ca(E,F)$, like, e.g., the characterization of
regular boundary points for the operator~$A$, which, actually, can be
expressed in terms of the capacity~$\cp$.

There are, however, problems where the capacity~$\ca$ can not be
replaced
by an equivalent capacity. An example is given by the study of the
asymptotic behaviour, as  ${j\to\infty}$, of the solutions~$u_j$ of the
Dirichlet problems
$$
\A{u_j}=f\,\ \hbox{in}\,\ \Om_j\,,\qquad
u_j=0\,\ \hbox{in}\,\ \partial\Om_j \,,
\leqno(\frm{hom})
$$
where~$A$ is a monotone operator of the form~(\rf{monop}),
$f\in\wm{\rn}$, and $(\Om_j)$ is a sequence of uniformly bounded
open sets in~$\rn$. Under some special assumptions on the structure of
the sets~$\Om_j$, this problem has been studied by means of the
capacities $\ca(E,F,s)$ in~[\rf{SKR-86}], [\rf{SKR-90}], [\rf{SKR-91}],
[\rf{SKR-93}], and~[\rf{SKR-94}], where a rigorous asymptotic
development of~$u_j$ is expressed in terms of the
$\ca$-potentials of suitable sets related to~$\Om_j$.

When~$A$ is the differential of a convex functional of the form
$$
G(u)=\int_{\rn} g(x,Du)\,dx \,,
$$
with $g(x,\cdot)$ even and homogeneous of degree~$p$, all assumptions
on the structure of the sets~$\Om_j$ can be avoided. In this case,
given a bounded open set~$\Om$ containing all sets~$\Om_j$, the
asymptotic behaviour of $(u_j)$ is determined by the limit, as
${j\to\infty}$, of the capacities $\ca({E\setminus\Om_j},\Om)$ on a
sufficiently large class of subsets~$E$ of~$\Om$ (see~[\rf{DM-DEF}]).
Since
these results depend strongly on the operator~$A$, it is clear that in
this problem~$\ca$ can not be replaced by an equivalent
capacity.

The properties of~$\ca$ will be used in a forthcoming paper to extend
the results of~[\rf{DM-DEF}] to the case of an arbitrary monotone
operator~$A$ of the form~(\rf{monop}). To this aim we intend to adapt
the techniques of~[\rf{DM-GAR-94}] to the non-linear case, by using the
results of the present paper and the compactness results
of~[\rf{DM-MUR-lum}],~[\rf{DM-MUR}], and~[\rf{CAS-GAR}].
\parag{\chp{notation}}{Notation and preliminaries}
\noindent
{\it Sobolev spaces and $p$-capacity.\/}
Let~$p$ and~$q$ be two real numbers with ${1<p\finite}$, ${1<q\finite}$,
and ${1/p+1/q=1}$.
For every open set~$\Om\subset\rn$ the Sobolev space $\wp{\Om}$
is defined as the space of all functions~$u$ in $\lp{\Om}$ whose first
order distribution derivatives $D_iu$ belong to $\lp{\Om}$, endowed
with the
norm
$$
\|u\|_{\wp{\Om}}^p\Eq \int_\Om |Du|^pdx \,+\, \int_\Om |u|^pdx \,,
$$
where $Du=(D_1u,\ldots,D_nu)$ is the gradient of~$u$. The space
$\wo{\Om}$ is the closure of $\cinftyo{\Om}$ in $\wp{\Om}$,
and $\wm{\Om}$ is the dual of $\wo{\Om}$. We shall always identify
each
function~$u$ of $\wo{\Om}$ with the function~$v$ of $\wp{\rn}$ such
that $v=u$ in~$\Om$ and $v=0$ in $\Om^c={\rn\setminus \Om}$. With
this
convention the space $\wo{\Om}$ can be regarded as a closed subspace
of
$\wp{\rn}$.

The lattice operations $\lor$ and $\land$ are defined by
$a\lor b=\max\{a,b\}$ and $a\land b=\min\{a,b\}$. It is well known
that, if~$u$ and~$v$ belong to $\wp{\Om}$, then ${u\lor v}$ and
${u\land v}$ belong to $\wp{\Om}$, and that the same property holds
for $\wo{\Om}$. The positive and the negative part of a function~$u$
are denoted by~$u^+$ and~$u^-$.

If $\Om\subset\rn$ is a bounded open set and $E\subset\Om$ is an
arbitrary set, the $p$-capacity of~$E$ in~$\Om$, denoted by
$\cp(E,\Om)$,  is defined as the infimum of $\int_\Om |Du|^pdx$ over
the set of all functions~$u$ in $\wo{\Om}$ such that $u\ge 1$ in a
neighbourhood of~$E$, with the usual convention
$\inf\emptyset=\pinfty$. It follows immediately from the definition
that
$$
\cp(E,\Om)\Eq \inf\{\cp(U,\Om): U\ \hbox{open,\ }
E \subset U\subset\Om \} \,.
\leqno(\frm{open})
$$
Moreover it is possible to prove that the set function $\cp(\cdot,\Om)$
is increasing and countably subadditive.

We say that a set~$N$ in~$\rn$ is {\it $\cp$-null\/} if
$\cp(N\cap\Om,\Om)=0$  for every bounded open set
$\Om\subset\rn$. It
is easy to prove that,  if~$N$ is contained in a bounded open
set~$\Om_0$
and  $\cp(N,\Om_0)=0$, then~$N$ is $\cp$-null, i.e.,
$\cp(N\cap\Om,\Om)=0$ for
every other bounded open set~$\Om$.

We say that a property ${\cal P}(x)$ holds {\it $\cp$-quasi
everywhere\/} (abbreviated as {\it $\cp$-q.e.\/}) in a set
$E\subset\rn$ if it
holds for all $x\in E$ except for a $\cp$-null set $N\subset E$. The
expression {\it almost everywhere\/} (abbreviated as {\it a.e.\/}) refers,
as usual, to the Lebesgue measure.

A function $u\colon\rn\to\r$ is said to be {\it quasi continuous\/} if
for every bounded open set $\Om\subset\rn$ and for every $\e>0$
there exists a set $E\subset\Om$, with $\cp(E,\Om)<\e$, such that the
restriction of~$u$ to ${\Om\setminus E}$ is continuous.

It is well known that every $u\in\wp{\rn}$ has a quasi continuous
representative, which is uniquely defined up to a $\cp$-null set.
In the sequel we shall always identify~$u$ with its  quasi continuous
representative, so that the pointwise values of a function
$u\in\wp{\rn}$ are  defined $\cp$-quasi everywhere in~$\rn$. We
recall that, if a sequence $(u_j)$ converges to~$u$ in $\wp{\rn}$, then a
subsequence of $(u_j)$ converges to~$u$ \qe~$\rn$. For all these
properties of quasi continuous representatives of Sobolev functions we
refer to~[\rf{EVA-GAR}], Section~4.8, [\rf{HEI-KIL-MAR}], Section~4,
[\rf{MAZ}], Section 7.2.4, and~[\rf{ZIE}], Section~3.

Given any set $E\subset\rn$ we define the Sobolev space $\wo{E}$ as
the set of all functions $u\in\wp{\rn}$ such that $u=0$ \qe~$E^c$, where
$E^c$ denotes the complement of~$E$ with respect to~$\rn$. It is easy to
see that $\wo{E}$ is a closed subspace of $\wp{\rn}$. The space
$\wm{E}$ is defined as the dual of $\wo{E}$ and the duality pairing is
denoted by $\langle\cdot,\cdot\rangle$. The transpose of the
imbedding of $\wo{E}$ into $\wp{\rn}$ defines a natural projection of
$\wm{\rn}$ onto $\wm{E}$, so that all elements of $\wm{\rn}$ can be
regarded as elements of $\wm{E}$. When~$E$ is open, our definitions
coincide with the classical definitions considered at the beginning of this
section (see~[\rf{HEI-KIL-MAR}], Theorem~4.5).

If~$f$ and~$g$ belong to $\wm{E}$, we say that $f=g$ in $\wm{E}$ if
$\langle f, v\rangle =\langle g, v\rangle$ for every $v\in\wo{E}$. We
say that $f\le g$ in $\wm{E}$ if
$\langle f, v\rangle \le \langle g, v\rangle$ for every $v\in\wo{E}$ with
$v\ge0$ \qe~$E$. It is easy to see that $f=g$ in $\wm{E}$ if and only if
$f\le g$  in $\wm{E}$ and $g\le f$ in $\wm{E}$.

The previous definitions allow us to consider the capacity $\cp(E,F)$
when~$E$ and~$F$ are arbitrary bounded sets in~$\rn$. In this case we
define
$$
\cp(E,F)=\min\,\{\int_F |Du|^pdx\,:\, u\in\wo{F},
\ u\ge1\,\ \cp\hbox{-q.e.\ in\ }\, E\}\,,
\leqno(\frm{cpmin})
$$
with the usual convention $\min\emptyset=\pinfty$. When~$F$ is open
and $E\subset F$, this definition agrees with the definition considered
at the beginning of the paper (see~[\rf{FED-ZIE}], Section~10, or
[\rf{HEI-KIL-MAR}], Corollary~4.13). If
$\cp(E,F)\finite$, then the minimum problem~(\rf{cpmin}) has a unique
minimum point, which is called the $\cp$-potential of~$E$ in~$F$.
\goodbreak\bigskip
\noindent
{\it Quasi open and quasi closed sets.\/}
We say that a set~$U$ in~$\rn$ is {\it $\cp$-quasi open\/}
(resp.\ {\it $\cp$-quasi closed\/}) if for every bounded open set
$\Om\subset\rn$ and for every $\e>0$ there exists an open
(resp.\ closed) set $V\subset\rn$ such that
$\cp({(U \bigtriangleup V)\cap\Om},\Om)<\e$, where $\bigtriangleup$
denotes the
symmetric difference of sets.

If a function $u\colon\rn\to\r$ is $\cp$-quasi continuous, then for
every $t\in\r$ the
sets $\{{u<t}\}=\{{x\in\rn}:{u(x)<t}\}$ and
$\{{u>t}\}=\{{x\in\rn}:{u(x)>t}\}$ are $\cp$-quasi open, while the sets
$\{{u\le t}\}=\{{x\in\rn}:{u(x)\le t}\}$ and
$\{{u\ge t}\}=\{{x\in\rn}:{u(x)\ge t}\}$ are $\cp$-quasi closed. In
particular this property holds for every $u\in\wp{\rn}$.

We shall frequently use the following lemma on the approximation of
the characteristic function of a $\cp$-quasi open set. We recall that
the characteristic function $1^{}_E$ of a set $E\subset\rn$ is defined
by $1^{}_E(x)=1$, if $x\in E$, and by $1^{}_E(x)=0$, if $x\in E^c$.
\lemma{\thm{approx}}{For every $\cp$-quasi open set~$U$ in~$\rn$
there exists an increasing sequence $(v_j)$ of
non-negative functions of $\wp{\rn}$ which converges to $1^{}_U$
$\cp$-quasi everywhere in~$\rn$.
}
\proof See~[\rf{DM-83}], Lemma~1.5, or~[\rf{DM-GAR}], Lemma~2.1.
\endproof
\spazio
The following lemma is used in the proof of Theorems~\rf{decreasingE}
and~\rf{increasingF}.
\lemma{\thm{approx2}}{Let~$U$ be the union of an increasing sequence
$(U_j)$ of $\cp$-quasi open sets in~$\rn$ and let~$F$ be an arbitrary set
in~$\rn$. Then for every function~$u$ of $\wo{F\cap U}$ there exists a
sequence $(u_j)$ which converges to~$u$ strongly in $\wp{\rn}$ and
such that $u_j\in\wo{F\cap U_j}$ for every~$j$.
}
\proof Let~$u$ be a function of $\wo{F\cap U}$. It is not restrictive to
assume that ${u\ge0}$ \qe~$\rn$. By Lemma~1.6 of~[\rf{DM-83}] there
exists an increasing sequence $(u_j)$ which converges to~$u$ strongly in
$\wp{\rn}$ and such that $u_j\le u\,1^{}_{U_j}$ \qe~$\rn$ for every~$j$.
Then the sequence $(u_j^+)$ converges to~$u$ strongly in $\wp{\rn}$
and satisfies $0\le u_j^+\le u\,1^{}_{U_j}$ \qe~$\rn$ for every~$j$. Since
$u=0$ \qe~$F^c$, we conclude that $u_j^+=0$ \qe{} ${F^c\cup U_j^c}$,
hence $u_j^+\in\wo{F\cap U_j}$ for every~$j$.
\endproof
\spazio
The following lemmas show that all $\cp$-quasi open sets and all
$\cp$-quasi closed sets can be approximated by an increasing sequence
of compact sets.
\lemma{\thm{compact1}}{Let~$U$ be a $\cp$-quasi open set in~$\rn$.
Then there exists an increasing sequence $(K_j)$ of compact subsets
of~$U$ whose union covers $\cp$-quasi all of~$U$.
}
\proof Since every $\cp$-quasi open set is the union of an increasing
sequence of $\cp$-quasi open bounded sets, we may assume that~$U$ is
bounded. Let~$\Om$ be a bounded open set in~$\rn$ containing~$U$.
Since~$U$ is $\cp$-quasi open, for every $k\in\n$ there exists an open
set~$V_k$, contained in~$\Om$, such that $\cp({U\bigtriangleup
V_k},\Om)<1/k$. By~(\rf{open}) there exists an open
set~$W_k$
such that ${U\bigtriangleup V_k}\subset
W_k\subset\Om$ and $\cp(W_k,\Om)<1/k$. This implies, in
particular, that $V_k\setminus W_k = U\setminus W_k$. As~$V_k$ is
open, it is the union of an increasing sequence $(C^j_k)_j$ of compact
sets. Let us define
$$
K_j=(C^j_1\setminus W_1) \cup \cdots \cup (C^j_j\setminus W_j) \,.
$$
Then $K_j$ is compact and contained in~$U$. As $C^j_k\subset
C^{j+1}_k$, the sequence $(K_j)$ is increasing. Moreover the union $E$ of
$(K_j)$ contains ${V_k\setminus W_k}={U\setminus W_k}$ for
every~$k$. Therefore
$\cp({U\setminus E},\Om)\le \cp(W_k,\Om) <1/k$  for every~$k$, and
hence $\cp({U\setminus E},\Om)=0$.
\endproof
\lemma{\thm{compact2}}{Let~$F$ be a $\cp$-quasi closed set in~$\rn$.
Then there exists an increasing sequence $(K_j)$ of compact subsets
of~$F$ whose union covers $\cp$-quasi all of~$F$.
}
\proof Since every $\cp$-quasi closed set is the union of an increasing
sequence of $\cp$-quasi closed bounded sets, we may assume that~$F$
is bounded. Let~$\Om$ be a bounded open set in~$\rn$
containing~$\overline F$. Since~$F$ is $\cp$-quasi closed, for every
$j\in\n$ there exists a compact set~$F_j$, contained in~$\Om$, such that
$\cp({F\bigtriangleup F_j},\Om)<2^{{-}j}$. By~(\rf{open}) there exists an
open set~$U_j$
such that ${F\bigtriangleup F_j}\subset
U_j\subset\Om$ and $\cp(U_j,\Om)<2^{{-}j}$. Let
$V_j=U_j\cup U_{j+1}\cup\cdots$, so that $V_{j+1}\subset V_j$ and
$\cp(V_j,\Om)<2^{1-j}$. Let $K_j$ be the compact set defined by
$K_j=F_j\setminus V_j$. As ${F\bigtriangleup F_j}\subset
U_j\subset V_j$, we have $K_j=F\setminus V_j$. This implies that
$K_j\subset F$ and that the sequence $(K_j)$ is increasing. Moreover the
union~$E$ of $(K_j)$ contains ${F\setminus V_j}$ for
every~$j$. Therefore
$\cp({F\setminus E},\Om)\le \cp(V_j,\Om) <2^{1-j}$  for every~$j$, and
hence $\cp({F\setminus E},\Om)=0$.
\endproof
\lemma{\thm{compact}}{Let~$E=U\cap F$ be the intersection of a
$\cp$-quasi open set~$U$ and a
$\cp$-quasi closed set~$F$.
Then there exists an increasing sequence $(K_j)$ of compact subsets
of~$E$ whose union covers $\cp$-quasi all of~$E$.
}
\proof The assertion follows from Lemmas~\rf{compact1}
and~\rf{compact2}.
\endproof
\goodbreak\bigskip
\noindent
{\it Measures.\/}
By a {\it Radon measure\/} on $\rn$ we
mean a continuous linear functional on the space $C_0(\rn)$
of all continuous functions with compact support in~$\rn$. It is
well known that for every Radon measure~$\lambda$ there exists a
countably additive set function~$\mu$, defined on the family of all
bounded Borel subsets of~$\rn$, such that
$\lambda(u)=\into u\,d\mu$ for every $u\in C_0(\Om)$. We shall
always identify the functional~$\lambda$ with the set
function~$\mu$.

We say that a Radon measure~$\mu$ on $\rn$ belongs to
$\wm{\rn}$ if there exists  $f\in\wm{\rn}$ such that
$$
\langle f,\varphi\rangle =\int_{\rn} \varphi\,d\mu
\qquad \forall \varphi\in\cinftyo{\rn}\,,
\leqno(\frm{measwm})
$$
where $\langle\cdot,\cdot\rangle$ denotes the duality pairing
between $\wm{\rn}$ and~$\wp{\rn}$. We shall always identify~$f$
and~$\mu$. Note that, by the Riesz Representation Theorem, for
every
non-negative functional $f\in\wm{\rn}$ there exists a
non-negative Radon measure~$\mu$ on~$\rn$ such
that~(\rf{measwm})
holds.

We say that a Radon measure~$\mu$ on~$\rn$ is
{\it $\cp$-absolutely continuous\/} if $\mu(N)=0$ for every
$\cp$-null Borel set $N\subset\rn$. It is well known that every
non-negative Radon measure~$\mu$ which belongs to $\wm{\rn}$
is
$\cp$-absolutely continuous and that, in this case, $\wp{\rn}\subset
L^1_\mu(\rn)$ and~(\rf{measwm}) holds for every
$\varphi\in\wp{\rn}$.

\goodbreak\bigskip
\noindent
{\it The monotone operator.\/} Let $a\colon \rn{\times}\rn\to \rn$ be a
function which satisfies the usual Carath\'eodory conditions, i.e., for
every $\xi\in\rn$ the function $x\mapsto \a{\xi}$ is (Lebesgue)
measurable on~$\rn$, and for a.e.\ $x\in\rn$ the function $\xi\mapsto
\a{\xi}$ is continuous on~$\rn$. We assume that there exist two
constants $c_1>0$ and $c_2>0$, and two functions $b_1\in\lone{\rn}$
and $b_2\in\lq{\rn}$, such that
$$
\ldisplaylinesno{\a{0}=0\,,
&(\frm{zero})
\cr
\aaa{\xi}{\eta}{\xi-\eta}\ge0\,,
&(\frm{mono})
\cr
\aa{\xi}{\xi}\ge c_1|\xi|^p - b_1(x) \,,
&(\frm{coerc})
\cr
|\a{\xi}| \le c_2|\xi|^{p-1} + b_2(x)
&(\frm{bound})
\cr}
$$
for a.e.\ $x\in\rn$ and for every $\xi$, $\eta\in\rn$. Note
that~(\rf{zero}) and~(\rf{mono}) imply
$$
\aa{\xi}{\xi}\ge 0
\leqno(\frm{positive})
$$
for a.e.\ $x\in\rn$ and for every $\xi\in\rn$.
Let
$A\colon\wp{\rn}\to\wm{\rn}$ be the operator defined by
$\A{u}=-\div\bigl(\a{Du}\bigr)$, i.e.,
$$
\langle\A{u},v\rangle \Eq \iaa{\rn}{Du}{Dv}
\leqno(\frm{operator})
$$
for every $u$, $v\in\wp{\rn}$.
\rem{\thm{properties}}{Since~$a$ satisfies the Carath\'eodory
conditions, by~(\rf{bound}) the operator~$A$ is continuous.
Moreover, for every $E\subset\rn$ and for every
$u\in\wp{\rn}$ we have
$$
\langle\A{u},v\rangle \Le  \bigl( c_2 \|Du\|_{\lpn{E}}^{p-1} \,+\,
\|b_2\|_{\lq{E}}\bigr) \|Dv\|_{\lpn{E}}
\qquad \forall v\in\wo{E} \,.
\leqno(\frm{bound1})
$$
 From~(\rf{mono}) we get
$$
\langle\A{u}-\A{v},u-v\rangle\ge0
\qquad\forall u\,,\ v\in\wp{\rn}\,,
\leqno(\frm{mono1})
$$
hence~$A$ is monotone. Inequality~(\rf{positive}) implies that
$$
\langle\A{u},u\rangle \ge 0
\qquad\forall u\in\wp{\rn}\,.
\leqno(\frm{positive1})
$$
By~(\rf{coerc}) for every $E\subset\rn$ we
have
$$
\langle\A{u},u\rangle \Ge c_1\int_{E}|Du|^pdx \,-\, \int_{E} b_1\,dx
\qquad\forall u\in\wo{E}\,.
\leqno(\frm{coerc1})
$$
By Poincar\'e's Inequality this implies that~$A$ is coercive on all
subspaces of the form $\{{u\in\wp{\rn}}:{u-\psi\in\wo{E}}\}$, with~$E$
bounded and $\psi\in\wp{\rn}$.
}

\parag{\chp{estimates}}{Some properties of the solutions}
In this section we prove some properties of the solutions~$u$ of the
Dirichlet problem
$$
\dirpb{u-\psi\in\wo{E}}
{\A{u}=f}
{\wm{E}},
\leqno(\frm{dir20})
$$
where $E$ is an arbitrary bounded set in $\rn$, $\psi$ is a function in
$\wp{\rn}$, and $f$ belongs to $\wm{E}$.
\th{\thm{existence}}{Let~$E$ be a bounded set in~$\rn$, let
$\psi\in\wp{\rn}$, and let $f\in\wm{E}$. Then the Dirichlet
problem~(\rf{dir20}) has at least a solution, and the set of all solutions
of~(\rf{dir20}) is bounded, closed, and convex in $\wp{\rn}$.
}
\proof By Remark~\rf{properties} the operator
$A\colon\wp{\rn}\to\wm{\rn}$ is
continuous and monotone.
Since~$E$ is bounded, $A$ is coercive on the set
$\{{u\in\wp{\rn}}:{u-\psi\in\wo{E}}\}$. Therefore the properties of the
set of the
solutions of~(\rf{dir20}) follow from the classical theory of monotone
operators (see, e.g.,~[\rf{KIN-STA}], Chapter~III).
\endproof
\lemma{\thm{lattice}}{Let~$E$ be a bounded set in~$\rn$, let
$\psi\in\wp{\rn}$, let $f\in\wm{E}$, and let~$u_1$ and~$u_2$ be two
solutions of~(\rf{dir20}). Then $u_1\lor u_2$ and $u_1\land u_2$ are
solutions of~(\rf{dir20}).
}
\proof Since $u_1-u_2\in\wo{E}$, by~(\rf{dir20}) we have
$$
\iaaa{E}{Du_1}{Du_2}{Du_1-Du_2} \Eq 0 \,.
$$
By~(\rf{mono}) we have $\aaa{Du_1}{Du_2}{Du_1-Du_2}\ge 0$ \ae~$E$,
hence
$$
\aaa{Du_1}{Du_2}{Du_1-Du_2}=0 \qquad \hbox{a.e.\ in\ }\, E\,.
\leqno(\frm0)
$$
Let us fix $v\in\wo{E}$ with $v\ge0$ \qe{}~$E$. For every $\e>0$ let us
define $v_\e=(\e v)\land({u_1-u_2})$. As $v_\e\in\wo{E}$,
by~(\rf{dir20}) we have
$$
\displaylines{
\iaa E{Du_1}{Dv_\e} \Eq \langle f,v_\e \rangle \,,
\cr
\iaa E{Du_2}{\e Dv - Dv_\e} \Eq \langle f,\e v - v_\e \rangle \,.
\cr}
$$
By adding these equalities we obtain
$$
\iaa E{Du_1}{Dv_\e} \,+\, \iaa E{Du_2}{\e Dv - Dv_\e} \Eq
\e \langle f,v \rangle \,.
$$
This implies that
$$
\displaylines{
\e \iaa{\{\e v<u_1-u_2\}}{Du_1}{Dv} \,+\,
\iaa{\{u_1-u_2\le \e v\}}{Du_1}{Du_1-Du_2} \,+
\cr
+\, \e \iaa{\{u_1-u_2\le \e v \}}{Du_2}{Dv} \,-\,
\iaa{\{u_1-u_2\le \e v\}}{Du_2}{Du_1-Du_2}  \Eq
\e \langle f,v \rangle \,.
\cr}
$$
By~(\rf0) we have
$$
\iaa{\{\e v<u_1-u_2\}}{Du_1}{Dv} \,+\,
\iaa{\{u_1-u_2\le \e v \}}{Du_2}{Dv} \Eq
\langle f,v \rangle \,.
$$
Passing to the limit as $\e\to0$ we get
$$
\iaa{\{u_1>u_2\}}{Du_1}{Dv} \,+\,
\iaa{\{u_1\le u_2\}}{Du_2}{Dv} \Eq
\langle f,v \rangle \,,
$$
which implies
$$
\iaa E{D(u_1\lor u_2)}{Dv}\Eq
\langle f,v \rangle
$$
for every $v\in\wo{E}$ with $v\ge0$ \qe~$E$. This proves that
$\A{(u_1\lor u_2)}=f$ in $\wm{E}$.  The equality
$\A{(u_1\land u_2)}=f$ in $\wm{E}$ can be proved in a similar way.
\endproof
\spazio
We are now in a position to prove the existence of a maximal and a
minimal solution of~(\rf{dir20}).
\th{\thm{maximal}}{Let~$E$ be a bounded set in~$\rn$, let
$\psi\in\wp{\rn}$, and let $f\in\wm{E}$. Then there exist two
solutions~$u_1$ and~$u_2$ of problem~(\rf{dir20}) such that $u_1\le
u\le u_2$ \qe~$E$ for every other solution~$u$ of~(\rf{dir20}).
}
\proof Let~$K$ be the set of all solutions of~(\rf{dir20}). By
Theorem~\rf{existence} $K$ is non-empty, bounded, closed, and convex
in $\wp{\rn}$. Since this space is separable, there exists a sequence
$(v_j)$
in~$K$ which is dense in~$K$. For every $k\in\n$ let us define
$$
u_1^k=\inf_{1\le j\le k} v_j\,,
\qquad
u_2^k=\sup_{1\le j\le k} v_j\,,
\qquad
u_1=\inf_{1\le j} v_j\,,
\qquad
u_2=\sup_{1\le j} v_j\,.
$$
By Lemma~\rf{lattice} both~$u_1^k$ and~$u_2^k$ belong to~$K$,
therefore the sequences $(u_1^k)$ and $(u_2^k)$ are bounded in
$\wp{\rn}$ and converge to~$u_1$ and~$u_2$ weakly in $\wp{\rn}$.
As~$K$ is weakly closed in $\wp{\rn}$, we conclude that~$u_1$
and~$u_2$ belong to~$K$, i.e., they are solutions of~(\rf{dir20}). If~$u$
is  another solution of~(\rf{dir20}), by the density of $(v_j)$ there
exists a  subsequence  $(v_{j_k})$ which converges to~$u$ strongly in
$\wp{\rn}$  and \qe~$\rn$. Since $u_1^{j_k}\le v_{j_k}\le u_2^{j_k}$
\qe~$\rn$ for every~$k$, we  conclude that $u_1\le u\le u_2$
\qe~$\rn$.
\endproof
\defin{\thm{maxsol}}{The functions~$u_1$ and~$u_2$ introduced in the
previous theorem are called the {\it minimal\/} and the {\it maximal
solution\/} of problem~(\rf{dir20}).
}
\spazio
The following lemma will be fundamental in the proof of the
monotonicity properties of the capacity~$\ca$ associated with the
operator~$A$.
\lemma{\thm{fundamental}}{Let~$B$ and~$C$ be two sets in~$\rn$, and
let $w_1$ and $w_2$ be two functions in
$\wp{\rn}$ such that $w_1=w_2$ \qe~$B$ and $w_1\le w_2$ \qe~$C$.
Assume that ${\A{w_1}=\A{w_2}}$ in
$\wm{C\setminus B}$. Then $\A{w_1}\ge\A{w_2}$ in $\wm{C}$.
}
\proof Let us fix a function~$v$ in $\wo{C}$ with $v\ge0$ \qe~$C$. For
every $\e>0$ let us define $v_\e=(\e v)\land(w_2-w_1)^+$.
Since $v_\e$ belongs to $\wo{{C\setminus B}}$ and $\A{w_1}=\A{w_2}$
in $\wm{{C\setminus B}}$, we have $\langle\A{w_1},v_\e\rangle=
\langle\A{w_2},v_\e\rangle$. Therefore by the monotonicity
condition~(\rf{mono})
$$
\displaylines{
\e \langle\A{w_1},v\rangle \,-\, \e \langle\A{w_2},v\rangle \Eq
\langle\A{w_1}-\A{w_2}, \e v-v_\e \rangle \Eq
\cr
\Eq \e \iaaa{C\cap\{w_2-w_1\le \e v\}}{Dw_1}{Dw_2}{Dv} \,-
\cr
-\, \iaaa{C\cap\{w_2-w_1\le \e v\}}{Dw_1}{Dw_2}{Dw_2-Dw_1} \Ge
\cr
\Ge \e \iaaa{C\cap\{w_2-w_1\le \e v\}}{Dw_1}{Dw_2}{Dv} \,.
\cr}
$$
Dividing by $\e$ and passing to the limit as $\e\to0$ we obtain
$$
\langle\A{w_1},v\rangle \,-\, \langle\A{w_2},v\rangle \Ge
\iaaa{C\cap\{w_2\le w_1\}}{Dw_1}{Dw_2}{Dv}\,.
\leqno(\frm{star})
$$
Since $w_1\le w_2$ \qe~$C$, we have $Dw_1=Dw_2$ \ae{}
$C\cap\{{w_2\le w_1}\} = C\cap\{{w_2= w_1}\}$. Therefore the right
hand side of~(\rf{star}) is equal to~$0$ and, consequently,
$\langle\A{w_1},v\rangle \ge \langle\A{w_2},v\rangle$,
which concludes the proof.
\endproof
\spazio
If~$u_1$ and~$u_2$ are two solutions of problem~(\rf{dir20}), then
$\A{u_1}=\A{u_2}$ in $\wm{E}$. As $u_1=u_2=\psi$ \qe~$E^c$, we have
also $\A{u_1}=\A{u_2}$ in $\wm{E^c}$. The following theorem shows
that actually $\A{u_1}=\A{u_2}$ in $\wm{\rn}$.
\th{\thm{b}}{Let~$E$ be a bounded set in~$\rn$, let
$\psi\in\wp{\rn}$, and let $f\in\wm{E}$. Then there exists
$g\in\wm{\rn}$ such that $\A{u}=g$ in $\wm{\rn}$ for every
solution~$u$ of problem~(\rf{dir20}).
}
\proof Let~$u_1$ and~$u_2$ be the minimal and the maximal solution
of~(\rf{dir20}). If we apply Lemma~\rf{fundamental} with $B=E^c$ and
$C=\rn$, we obtain that $\A{u_1}\ge\A{u}\ge\A{u_2}$ in $\wm{\rn}$
for every solution~$u$ of~(\rf{dir20}). Therefore it is enough to prove
that $\A{u_1}\le\A{u_2}$ in $\wm{\rn}$. Let $\varphi\in\cinftyo{\rn}$
with $\varphi\ge0$ in~$\rn$ and let $c=\|\varphi\|^{}_{\linfty{\rn}}$.
Let us fix a function $\chi$ in $\cinftyo{\rn}$ such that $\chi=c$ in
$E\cup\{{\varphi\neq0}\}$ and $\chi\ge0$ in~$\rn$.
Then $\chi-\varphi$ belongs to
$\cinftyo{\rn}$ and $\chi-\varphi\ge0$ in~$\rn$. As
$\A{u_1}\ge\A{u_2}$ in  $\wm{\rn}$, we have
$$
\langle \A{u_1}, \chi-\varphi \rangle \Ge
\langle  \A{u_2}, \chi-\varphi \rangle \,.
\leqno(\frm{chi-phi})
$$
Since $u_1=u_2=\psi$ \qe~$E^c$, we have $Du_1=Du_2$ \ae{}
$\{{\chi\neq c}\}$. As $D\chi=0$ \ae{} $\{{\chi=c}\}$, we have
$$
\langle \A{u_1} ,\chi\rangle \Eq
\iaa{\{{\chi\neq c}\}}{Du_1}{D\chi} \Eq
\iaa{\{{\chi\neq c}\}}{Du_2}{D\chi} \Eq
\langle \A{u_2} ,\chi\rangle \,.
$$
By~(\rf{chi-phi}) this implies $ \langle \A{u_1} ,\varphi\rangle \le
\langle \A{u_2} ,\varphi\rangle$
for every $\varphi\in\cinftyo{\rn}$ with $\varphi\ge0$ and, by
density, this implies
$\A{u_1}\le\A{u_2}$ in $\wm{\rn}$.
\endproof
\cor{\thm{restriction}}{Let~$E$ be a bounded set in~$\rn$, let
$\psi\in\wp{\rn}$, let $f\in\wm{E}$, and let~$u_1$ and~$u_2$ be the
minimal and the maximal solution of~(\rf{dir20}). Let~$B$ be a subset
of~$E$. Then~$u_1$ coincides with the minimal solution~$v_1$ of the
Dirichlet problem
$$
\dirpb{v-u_1\in\wo{B}}{\A{v}=f}{\wm{B}},
\leqno(\frm{r1})
$$
and~$u_2$ coincides with the maximal solution~$v_2$ of the Dirichlet
problem
$$
\dirpb{v-u_2\in\wo{B}}{\A{v}=f}{\wm{B}}.
\leqno(\frm{r2})
$$
}
\negbigskip
\proof It is clear that~$u_1$ is a solution of~(\rf{r1}). Let~$v_0$ be
another solution of~(\rf{r1}). We have to prove that $u_1 \le v_0$
\qe~$B$. By Theorem~\rf{b}, applied to problem~(\rf{r1}), there exists
$g\in\wm{\rn}$ such that $\A{v}=g$ in $\wm{\rn}$ for every
solution~$v$ of~(\rf{r1}). In particular we have $\A{v_0}=\A{u_1}$ in
$\wm{\rn}$.  Since~$u_1$ is a solution of~(\rf{dir20}), we have
$\A{u_1}=f$ in $\wm{E}$. This implies that $\A{v_0}=f$ in $\wm{E}$,
thus~$v_0$ is a solution of~(\rf{dir20}) too. By the minimality of~$u_1$
we conclude that $u_1\le v_0$ \qe~$E$, and, therefore, $u_1$ is the
minimal solution of~(\rf{r1}). The proof for~$u_2$ is analogous.
\endproof
\cor{\thm{b1}}{Let~$E$ be a bounded set in~$\rn$, let
$\psi\in\wp{\rn}$, let $f\in\wm{\rn}$, and let~$u_1$ and~$u_2$ be two
solutions of~(\rf{dir20}). Then $\langle{\A{u_1}-f},u_1\rangle =
\langle{A{u_2}-f},u_2\rangle$.
}
\proof Since $u_1-u_2$ belongs to $\wo{E}$, by~(\rf{dir20}) we have
$\langle{\A{u_1}-f},{u_1-u_2}\rangle = 0$, hence
$\langle{\A{u_1}-f},u_1\rangle = \langle{\A{u_1}-f},u_2\rangle$. By
Theorem~\rf{b} we have $\A{u_1}=\A{u_2}$ in $\wm{\rn}$, hence
$\langle{\A{u_1}-f},u_2\rangle = \langle{\A{u_2}-f},u_2\rangle$.
\endproof
\spazio
The following lemma will be used in the proof of the Comparison
Principle.
\lemma{\thm{subsolution}}{Let~$E$ be a bounded open set in~$\rn$, let
$w$, $\psi\in\wp{\rn}$, and let $f\in\wm{E}$. Assume that
$\A{w}\le f$ in $\wm{E}$ and that $(w-\psi)^+\in\wo{E}$, i.e.,
$w\le\psi$ \qe~$E^c$. Then there exists a solution~$u$ of~(\rf{dir20})
such that $u\ge w$ \qe~$E$.
}
\proof Let $K$ be the set of all functions~$v$ in $\wp{\rn}$ such that
$v=\psi$ \qe~$E^c$ and $v\ge w$ \qe~$E$. As $({w-\psi})^+$ belongs
to $\wo{E}$, the function $\psi+({w-\psi})^+$ belongs to~$K$, so that~$K$
is non-empty. By Remark~\rf{properties}, using the classical theory of
monotone operators (see, e.g.,~[\rf{KIN-STA}], Chapter~III), we can find
a
solution~$u$ of the variational inequality
$$
\cases{u\in K\,,&
\cr\cr
\langle\A{u},v-u\rangle \ge \langle f,v-u\rangle
\qquad \forall v\in K\,.&
\cr}
\leqno(\frm{varin})
$$
We want to prove that~$u$ is a solution of~(\rf{dir20}). If $z$ belongs to
$\wo{E}$ and $z\ge0$ \qe~$E$, then $v=u+z$ belongs to~$K$, hence
$$
\langle\A{u},z\rangle \ge \langle f,z\rangle \,.
\leqno(\frm{pos})
$$
In order to prove the opposite inequality, for every $\e>0$ we consider
the function $z_\e=(\e z)\land(u-w)$. Since $u=\psi\ge w$ and $z=0$
\qe~$E^c$, we have $z_\e=0$ \qe~$E^c$, hence $z_\e\in\wo{E}$. As
$u-z_\e \ge u - ({u-w})=w$ \qe~$E$, the function $v_\e=u-z_\e$ belongs
to~$K$. Therefore (\rf{varin}) yields
$\langle\A{u},z_\e\rangle \le \langle f,z_\e\rangle$. This implies
$$
\e \iaa{\{\e z<u-w\}}{Du}{Dz} \,+\,
\iaa{\{\e z\ge u-w\}}{Du}{Du-Dw} \Le
\langle f,z_\e\rangle \,.
\leqno(\frm{ze})
$$
Since $\A{w}\le f$ in $\wm{E}$ and $\e z - z_\e\ge0$ \qe~$E$, we have
also $\langle\A{w},{\e z-z_\e}\rangle \le \langle f,{\e z-z_\e}\rangle$,
which gives
$$
\e \iaa{\{\e z\ge u-w\}}{Dw}{Dz} \,-\,
\iaa{\{\e z\ge u-w\}}{Dw}{Du-Dw} \Le
\langle f,{\e z-z_\e}\rangle \,.
$$
By adding this inequality to~(\rf{ze}) we obtain
$$
\displaylines{
\e \iaa{\{\e z<u-w\}}{Du}{Dz} \,+\,
\e \iaa{\{\e z\ge u-w\}}{Dw}{Dz} \,+
\cr
+\, \iaaa{\{\e z\ge u-w\}}{Du}{Dw}{Du-Dw} \Le
\e \langle f,z \rangle \,.
\cr}
$$
By the monotonicity condition~(\rf{mono}) we get
$$
\iaa{\{\e z<u-w\}}{Du}{Dz} \,+\,
\iaa{\{\e z\ge u-w\}}{Dw}{Dz} \Le
 \langle f,z \rangle \,,
$$
and taking the limit as $\e\to0$ we obtain
$$
\iaa{\{u>w\}}{Du}{Dz} \,+\,
\iaa{\{u\le w\}}{Dw}{Dz} \Le
 \langle f,z \rangle \,.
\leqno(\frm{zed})
$$
As $u\ge w$ \qe~$\rn$, we have $Du=Dw$ \ae{}
$\{{u\le w}\}= \{{u= w}\}$, so that~(\rf{zed}) gives
$\langle \A{u},z \rangle\le \langle f,z \rangle$. Together with~(\rf{pos})
this implies $\A{u}=f$ in $\wm{E}$. As $u\in K$, we have also
${u-\psi}\in\wo{E}$ and $u\ge w$ \qe~$E$ as required.
\endproof
\rem{\thm{supersolution}}{If in Lemma~\rf{subsolution} we assume
that $\A{w}\ge f$ in $\wm{E}$ and that $(w-\psi)^-\in\wo{E}$, i.e.,
$w\ge\psi$ \qe~$E^c$, then we can prove that there exists a
solution~$u$ of~(\rf{dir20}) such that $u\le w$ \qe~$E$.
}
\spazio
We are now in a position to prove the Comparison Principle.
\th{\thm{comparison}}{Let~$E$ be a bounded set in~$\rn$, let
$\varphi$, $\psi\in \wp{\rn}$, and let $f$, $g\in\wm{E}$. Assume that
$f\le g$ in $\wm{E}$ and $(\varphi-\psi)^+\in \wo{E}$, i.e.,
$\varphi\le\psi$ \qe~$E^c$. Let~$u$ and~$v$ be two
solutions of the Dirichlet problems
$$
\dirpb{u-\varphi\in\wo{E}}{\A{u}=f}{\wm{E}},
\qquad\qquad
\dirpb{v-\psi\in\wo{E}}{\A{v}=g}{\wm{E}},
\leqno(\frm{dir11})
$$
let~$u_1$ and~$v_1$ be the minimal
solutions, and let~$u_2$ and~$v_2$ be the maximal solutions. Then
$u_1\le v_1\le v$
and $u\le u_2\le v_2$ \qe~$E$.
}
\proof Since $\A{u_2}=f\le g$ in $\wm{E}$ and $u_2=\varphi\le\psi$
\qe~$E^c$, by Lemma~\rf{subsolution}  there exists a solution~$v_0$ of
the second problem in~(\rf{dir11}) such that $u_2\le v_0$ \qe~$E$.
Since~$u_2$ and~$v_2$ are the maximal solutions, we have $u\le u_2$
and $v_0\le v_2$ \qe~$E$, hence $u\le u_2\le v_2$ \qe~$E$. The other
inequalities can be proved in the same way by using
Remark~\rf{supersolution}.
\endproof
\parag{\chp{capacity}}{Capacity, capacitary potentials, and capacitary
distributions}
In this section we introduce the capacity~$\ca$ associated with the
operator~$A$, and the related notions of $\ca$-potential and
$\ca$-distributions.
\defin{\thm{compatible}}{We say that two bounded sets~$E$ and~$F$
are {\it $\cp$-compatible\/}
if there exists a function $\psi$ in $\wp{\rn}$ such that $\psi=1$
\qe~$E$ and $\psi=0$ \qe~$F^c$.
}
\rem{\thm{compatible1}}{It is clear that~$E$ and~$F$ are
$\cp$-compatible if and only if there exists a function $\psi$ in
$\wp{\rn}$ such that $0\le\psi\le1$ \qe~$\rn$, $\psi=1$ \qe~$E$, and
$\psi=0$ \qe~$F^c$.
}
\rem{\thm{compatible2}}{If $\overline E\subset \interior F$,
then~$E$ and~$F$ are
$\cp$-compatible. The converse is not true. For instance, for every
function~$u$ in $\wp{\rn}$ the sets $E=\{{u>t}\}$ and $F=\{{u>s}\}$ are
$\cp$-compatible for $s<t$, with $\psi=(t-s)^{{-}1}((u\land t)-s)^+$, but,
in general, $\overline E$ is not contained in $\interior F$.
}
\rem{\thm{compatible3}}{If~$E_1$ and~$E_2$ are
$\cp$-compatible with~$F$, so is ${E_1\cup E_2}$. For the proof, let
us consider two functions $\psi_1$ and $\psi_2$ as in
Remark~\rf{compatible1} with $E=E_1$ and $E=E_2$. Then
$\psi_1\lor\psi_2$ satisfies the same conditions with
$E={E_1\cup E_2}$. Similarly we can prove that, if~$E$ is
$\cp$-compatible with~$F_1$ and~$F_2$, then~$E$ is
$\cp$-compatible with ${F_1\cap F_2}$.
}
\defin{\thm{potential}}{Let~$E$ and~$F$ be two $\cp$-compatible
bounded sets in~$\rn$ and let $\psi$ be a function as in
Definition~\rf{compatible}. Every solution~$u$ of the Dirichlet problem
$$
\dirpb{u-\psi\in\wo{F\setminus E}}
{\A{u}=0}
{\wm{F\setminus E}},
\leqno(\frm{pot})
$$
is called a {\it $\ca$-potential\/} of~$E$ in~$F$. The maximal and the
minimal solutions of~(\rf{pot}) are called the {\it maximal\/} and {\it
minimal $\ca$-potentials\/} of~$E$ in~$F$.
}
\rem{\thm{potential1}}{Clearly the definition of $\ca$-potential does
not depend on the choice of~$\psi$. By the definition of the space
$\wo{F\setminus E}$ and by the properties of $\psi$ we have that
every $\ca$-potential~$u$ of~$E$ in~$F$ satisfies $u=1$ \qe~$E$ and
$u=0$ \qe~$F^c$.
}
\rem{\thm{potential2}}{Let~$u_1$ and~$u_2$ be the minimal and the
maximal $\ca$-potentials of~$E$ in~$F$. Since $\a{0}=0$ by~(\rf{zero}),
the Comparison Principle (Theorem~\rf{comparison})  implies that $0\le
u_2$ and $u_1\le1$ \qe{} ${F\setminus E}$.
}
\defin{\thm{cap}}{Let~$E$ and~$F$ be two
bounded sets in~$\rn$. If~$E$ and~$F$ are $\cp$-compatible, the
{\it capacity of~$E$ in~$F$ relative to the operator\/}~$\A{}$ is defined
as
$$
\ca(E,F)\Eq \langle \A{u},u\rangle \Eq
\iaa{F\setminus E}{Du}{Du}\,,
$$
where~$u$ is any $\ca$-potential of~$E$ in~$F$. By Corollary~\rf{b1}
this definition is independent of the choice of~$u$. If~$E$ and~$F$ are
not $\cp$-compatible, we define ${\ca(E,F)=\pinfty}$.
}
\rem{\thm{cap1}}{By~(\rf{positive1}) we have $\ca(E,F)\ge0$, and
by~(\rf{zero}) we have $\ca(\emptyset,F)=0$.
By~(\rf{pot}) we have also
$\ca(E,F)= \langle\A{u}, v\rangle$ for every $\ca$-potential~$u$ of~$E$
in~$F$ and for every function~$v$ in $\wp{\rn}$ such that $v=1$
\qe~$E$ and $v=0$ \qe~$F^c$.
}
\rem{\thm{cap3}}{It follows immediately from the definitions that if
$E_1$, $E_2$, $F_1$, $F_2$ are bounded sets in~$\rn$ and
$F_1\bigtriangleup F_2$ and $E_1\bigtriangleup E_2$ are
$\cp$-null sets, then $\ca(E_1,F_1)=\ca(E_2,F_2)$ and the
$\ca$-potentials are the same.
}
\spazio
If $\A{u}=-\div\bigl(|Du|^{p-2}Du\bigr)$, then $\ca=\cp$. In the general
case the relationship between $\ca$ and $\cp$ is given by the following
proposition.
\prop{\thm{ca-cp}}{Let~$E$ and~$F$ be two bounded sets in~$\rn$.
Then
$$
\ldisplaylinesno{
\ca(E,F) \Ge c_1\,\cp(E,F)\, -\, \|b_1\|_{\lone{F}} \,,
&(\frm{lower})
\cr
\ca(E,F) \Le k_1 \cp(E,F) \,+\, k_2(F) \,\cp(E,F)^{1/p} \,,
&(\frm{upper})
\cr}
$$
with
$$
k_1\Eq { (4 c_2)^p \over p(q c_1)^{p-1} }\,,
\qquad
k_2(F)= {4 c_2\over c_1^{1/q}}\|b_1\|_{\lone{F}}^{1/q}\, + \,
4\|b_2\|^{}_{\lq{F}}\,,
\leqno(\frm{const})
$$
where $c_1$, $c_2$ and $b_1$, $b_2$ are the constants and the
functions which
appear in~(\rf{coerc}) and~(\rf{bound}). If~$b_1$ and~$b_2$ belong to
$\linfty{F}$, then
$$
\ca(E,F) \Le \bigl(k_1 + k_3(F)\bigr) \,\cp(E,F) \,,
\leqno(\frm{upper2})
$$
with
$$
k_3(F) \Eq
2^{p+1} \Big({ c_2\over c_1^{1/q}}\|b_1\|_{\linfty{F}}^{1/q}\, + \,
\|b_2\|^{}_{\linfty{F}}\Bigr) {\rm diam}(F)^{p-1} \,,
\leqno(\frm{const2})
$$
where ${\rm diam}(F)$ is the diameter of the set~$F$.
}
\proof To prove~(\rf{lower}) we may assume that $\ca(E,F)\finite$.
Let~$u$ be a $\ca$-potential of~$E$ in~$F$. By~(\rf{coerc1}) we
have
$$
\ca(E,F) \Eq \langle \A{u},u \rangle \Ge
c_1 \int_F |Du|^pdx \,-\, \int_F b_1\,dx \,.
\leqno(\frm{a0})
$$
Since by~(\rf{cpmin})
$$
\cp(E,F) \Le \int_F |Du|^pdx\,,
$$
from~(\rf{a0}) we obtain~(\rf{lower}).

To prove~(\rf{upper}) and~(\rf{upper2}) we may assume that
$\cp(E,F)\finite$.
Let~$w$ be the $\cp$-potential of~$E$ in~$F$, let $v=(2w-1)^+$, and let
$G=\{{v>0}\}=\{{w>{1\over 2}}\}$. Since ${w=1}$ \qe~$E$ and ${w=0}$
\qe~$F^c$, we have ${v=1}$ \qe~$E$ and ${v=0}$ \qe~$F^c$.
 From~(\rf{bound1}) and from Remark~\rf{cap1} we obtain
$$
\ca(E,F) \Eq \langle\A{u}, v\rangle \Le
\bigl( c_2 \| Du \|_{\lpn{G}}^{p-1} + \|b_2\|^{}_{\lq{G}} \bigr)
\|Dv\|^{}_{\lpn{G}} \,.
\leqno(\frm{a1})
$$
By~(\rf{coerc1}) we have
$$
c_1^{1/q} \| Du \|_{\lpn{F}}^{p-1} \Le  \ca(E,F)^{1/q} \,+\,
\|b_1\|_{\lone{F}}^{1/q} \,,
\leqno(\frm{a2})
$$
while the definition of~$v$ gives
$$
\| Dv \|^{}_{\lpn{G}} \Eq 2\,\| Dw \|^{}_{\lpn{G}}  \Le 2\,\cp(E,F)^{1/p} \,.
\leqno(\frm{a3})
$$
 From~(\rf{a1}), (\rf{a2}), and~(\rf{a3}) we get
$$
\ca(E,F) \Le 2\Bigl( {c_2 \over c_1^{1/q}} \ca(E,F)^{1/q} +
{c_2 \over c_1^{1/q}} \|b_1\|_{\lone{G}}^{1/q} +
\|b_2\|^{}_{\lq{G}} \Bigr) \cp(E,F)^{1/p} \,.
$$
Using Young's Inequality 
we obtain
$$
\ca(E,F) \Le { (4 c_2)^p \over p(q c_1)^{p-1} } \cp(E,F) \,+\,
\Bigl( {4 c_2\over c_1^{1/q}}\|b_1\|_{\lone{G}}^{1/q}\, + \,
4\|b_2\|^{}_{\lq{G}} \Bigr) \cp(E,F)^{1/p} \,,
\leqno(\frm{long})
$$
which implies~(\rf{upper}). If $b_1$ and $b_2$ belong to $\linfty{F}$,
then
$$
\displayliness{
\|b_1\|_{\lone{G}}^{1/q} \le \meas(G)^{1/q} \|b_1\|_{\linfty{G}}^{1/q} \,,
\cr
\|b_2\|^{}_{\lq{G}} \le \meas(G)^{1/q} \|b_2\|^{}_{\linfty{G}} \,.
\cr}
\leqno(\frm{measG})
$$
As $G=\{{w>{1\over 2}}\}$, by Poincar\'e's Inequality we have
$$
\meas(G)\Le 2^p \|w\|_{\lp{F}}^p \Le
2^p{\rm diam}(F)^p\|Dw\|_{\lp{F}}^p \Eq 2^p{\rm diam}(F)^p \cp(E,F) \,.
$$
Therefore~(\rf{measG}) implies
$$
\displaylines{
\|b_1\|_{\lone{G}}^{1/q} \le
\|b_1\|_{\linfty{G}}^{1/q} 2^{p-1}{\rm diam}(F)^{p-1} \cp(E,F)^{1/q} \,,
\cr
\|b_2\|^{}_{\lq{G}} \le
\|b_2\|^{}_{\linfty{F}} 2^{p-1}{\rm diam}(F)^{p-1} \cp(E,F)^{1/q} \,,
\cr}
$$
which, together with~(\rf{long}), yields~(\rf{upper2}).
\endproof
\rem{\thm{cap2}}{If $b_1=0$ \ae~$F$ and $b_2\in \linfty{F}$, then for
every bounded set~$E$ we have
$$
c_1\,\cp(E,F) \Le \ca(E,F) \Le \bigl(k_1 + k_3(F)\bigr) \, \cp(E,F) \,,
$$
where~$c_1$, $k_1$, and~$k_3(F)$  are defined in~(\rf{coerc}),
(\rf{const}), and~(\rf{const2}).
}
\spazio
The following lemma will be useful in the proof of
Theorems~\rf{capdistr1}, \rf{increasingE}, \rf{decreasingE},
and~\rf{caps15}.
\lemma{\thm{signA}}{Let~$E$ and~$F$ be two
$\cp$-compatible bounded sets in~$\rn$ and let~$u$ be a
$\ca$-potential of~$E$ in~$F$. Then $\A{u}\ge0$ in $\wm{F}$ and
$\A{u}\le0$ in $\wm{E^c}$.
}
\proof Let~$u_1$ and~$u_2$ be the minimal and the maximal
$\ca$-potentials of~$E$ in~$F$. By Remark~\rf{potential2} we have
$u_1\le 1$
\qe~$\rn$. If we apply Lemma~\rf{fundamental} with
$B=E$, $C=F$, $w_1=u_1$, and we choose as $w_2$ any function in
$\wp{\rn}$ which is equal to~$1$ \qe~$F$, we obtain
$\A{u_1}\ge\A{w_2}=0$ in $\wm{F}$. Since $\A{u}=\A{u_1}$ in
$\wm{\rn}$ by Theorem~\rf{b}, we conclude that $\A{u}\ge0$ in
$\wm{F}$.

By Remark~\rf{potential2} we have
$u_2\ge 0$ \qe~$\rn$. If we apply Lemma~\rf{fundamental}
with $B=F^c$, $C=E^c$, $w_1=0$, and $w_2=u_2$, we obtain
$\A{u_2}\le\A{0}=0$ in $\wm{E^c}$. Since $\A{u}=\A{u_2}$ in
$\wm{\rn}$ by Theorem~\rf{b}, we conclude that $\A{u}\le0$ in
$\wm{E^c}$.
\endproof
\spazio
By Theorem~\rf{b} there exists
$g\in\wm{\rn}$ such that $\A{u}=g$ in $\wm{\rn}$ for every
$\ca$-potential~$u$ of~$E$ in~$F$. The following
theorem gives a precise description of~$g$.
\th{\thm{capdistr1}}{Let~$E$ and~$F$ be two $\cp$-compatible
bounded sets in~$\rn$. Then there exists a unique pair $(\lambda,\nu)$
of non-negative Radon measures on~$\rn$ such that:
\smallskip
\item{\rm(a)} $\lambda$ and~$\nu$ are mutually singular;
\smallskip
\item{\rm(b)} for every $\ca$-potential~$u$ of~$E$ in~$F$ and for every
$\varphi\in\cinftyo{\rn}$ we have
$$
\langle \A{u}, \varphi\rangle \Eq
\int_{\rn}\varphi\,d\lambda \,-\, \int_{\rn}\varphi\,d\nu \,.
$$
\smallskip\noindent
Moreover, the following properties hold:
\smallskip
\item{\rm(c)} the measures~$\lambda$ and~$\nu$ are bounded and
$\cp$-absolutely continuous;
\smallskip
\item{\rm(d)} condition~{\rm(b)} holds also for every
$\varphi\in\wp{\rn}\cap\linfty{\rn}$;
\smallskip
\item{\rm(e)} the measure $\lambda-\nu$ belongs to $\wm{\rn}$;
\smallskip
\item{\rm(f)} $\supp\,\lambda\subset\partial E$ and
$\supp\,\nu\subset\partial F$;
\smallskip
\item{\rm(g)}
$\displaystyle v\in\wo{E}\cup\wo{E^c} \ \Longrightarrow \ v=0$
$\lambda$-\ae~$\rn$;
\smallskip
\item{\rm(h)}
$\displaystyle v\in\wo{F}\cup\wo{F^c}\ \Longrightarrow \ v=0$
$\nu$-\ae~$\rn$;
\smallskip
\item{\rm(i)} $\lambda(U)=0$ whenever
$U\cap E$ and $U\cap E^c$ are $\cp$-quasi open;
\smallskip
\item{\rm(j)} $\nu(U)=0$ whenever $U\cap F$ and $U\cap F^c$ are
$\cp$-quasi open;
\smallskip\nobreak
\item{\rm(k)} $\lambda(F^c)=0$ and $\nu(E)=0$.
\par
}
\proof By Theorem~\rf{b} there exists $g\in\wm{\rn}$ such that
$\A{u}=g$ in $\wm{\rn}$ for every $\ca$-potential~$u$ of~$E$ in~$F$.
We want to prove that there exists a unique pair $(\lambda,\nu)$
of mutually singular non-negative Radon measures on~$\rn$ such that
$$
\langle g, \varphi\rangle \Eq
\int_{\rn}\varphi\,d\lambda \,-\, \int_{\rn}\varphi\,d\nu
\qquad \forall \varphi\in\cinftyo{\rn} \,,
\leqno(\frm{cinftyo})
$$
and that~$\lambda$ and~$\nu$ satisfy properties (c)--(k).

Let us fix a function $\psi$ as in Remark~\rf{compatible1}. By
Lemma~\rf{signA} we have $g\ge0$ in $\wm{F}$. Since $\psi\varphi$
belongs to $\wo{F}$ for every $\varphi$ in $\cinftyo{\rn}$, we conclude
that
$$
\langle g,\psi\varphi\rangle\ge0
\qquad\forall \varphi\in\cinftyo{\rn}\,,\ \varphi\ge0\,.
$$
By the Riesz Representation Theorem there exists a non-negative
Radon measure~$\lambda$ on~$\rn$ such that
$$
\langle g ,\psi\varphi\rangle \Eq \int_{\rn} \varphi\,d\lambda
\qquad\forall \varphi\in\cinftyo{\rn}\,.
\leqno(\frm{nu})
$$

In order to construct $\nu$, we recall that by
Lemma~\rf{signA} we have $g\le0$ in $\wm{E^c}$.
Since $(1-\psi)\varphi$ belongs to
$\wo{E^c}$ for every $\varphi$ in $\cinftyo{\rn}$, we conclude that
$$
\langle g,(1- \psi)\varphi\rangle\le0
\qquad\forall \varphi\in\cinftyo{\rn}\,,\ \varphi\ge0\,.
$$
By the Riesz Representation Theorem there exists a non-negative
Radon measure~$\nu$ on~$\rn$ such that
$$
\langle g ,(1-\psi)\varphi\rangle \Eq -\int_{\rn} \varphi\,d\nu
\qquad\forall \varphi\in\cinftyo{\rn}\,.
\leqno(\frm{lambda})
$$
 From~(\rf{nu}) and~(\rf{lambda}) we obtain
$$
\langle g,\varphi\rangle\Eq
\langle g,\psi\varphi\rangle \,+\,
\langle g,(1-\psi)\varphi\rangle
\Eq \int_{\rn}\varphi\,d\lambda \,-\,
\int_{\rn}\varphi\,d\nu \,,
$$
which proves~(\rf{cinftyo}) and hence~(b). Property~(e) follows from~(b)
and from the fact that $\A{u}$ belongs to $\wm{\rn}$.

Let us prove~(c). As ${g=0}$ in $\wm{F^c}$ by~(\rf{zero}),
using~(\rf{nu}) and~(\rf{lambda}) we obtain
$$
\int_{\rn} \varphi\,d\lambda \Eq
\int_{\rn} \varphi\,d\nu \Eq 0
$$
for every $\varphi\in\cinftyo{\rn}$ with
$\supp\,\varphi \subset F^c$.
Therefore the supports of~$\lambda$ and~$\nu$ are contained in the
compact set $\overline F$. This implies that the measures~$\lambda$
and~$\nu$ are bounded. It remains to show that~$\lambda$
and~$\nu$ vanish on all $\cp$-null sets. To this aim, it is sufficient
to prove that $\lambda(C)=\nu(C)=0$ for every $\cp$-null compact set
$C\subset\rn$. In this case it is possible to construct a
sequence $(\varphi_j)$ of functions in $\cinftyo{\rn}$ such that
$0\le\varphi_j\le1$ in~$\rn$, $\varphi_j=1$ in~$C$, and $(\varphi_j)$
converges to~$0$ strongly in $\wp{\rn}$. Then by~(\rf{nu})
 for every~$j$ we have
$$
\lambda(C)\Le \int_{\rn} \varphi_j\,d\lambda \Eq
\langle g,\psi\varphi_j\rangle \,.
$$
Since $(\psi\varphi_j)$ converges to~$0$ strongly in $\wp{\rn}$,
passing
to the limit as $j\to\infty$ we obtain $\lambda(C)=0$. In the same way
we
prove that $\nu(C)=0$.

Since the measures~$\lambda$ and~$\nu$ are bounded and
$\cp$-absolutely continuous, every function of
$\wp{\rn}\cap\linfty{\rn}$ belongs to $L^1_\lambda(\rn)$ and to
$L^1_\nu(\rn)$, and thus, by an easy approximation argument,
from~(\rf{nu}) and~(\rf{lambda}) we obtain
$$
\langle g,\psi v\rangle \Eq \int_{\rn}v\,d\lambda\,,
\quad \qquad
\langle g,(1-\psi) v\rangle \Eq -\int_{\rn}v\,d\nu
\leqno(\frm{psi})
$$
for every $v\in\wp{\rn} \cap\linfty{\rn}$, which implies~(d).

By considering separately the positive and the negative part of~$v$, it is
enough to prove~(g) when~$v$ is non-negative. Let us fix
$v\in\wo{E}\cup\wo{E^c}$ with $v\ge 0$ \qe~$\rn$.
If~$v$ belongs to $\wo{E}$, then $v=\psi v$ \qe~$\rn$. Since ${g=0}$ in
$\wm{E}$ by~(\rf{zero}), from~(\rf{psi}) we obtain
$$
\int_{\rn}v\,d\lambda \Eq \langle g, v\rangle \Eq  0 \,.
$$
If~$v$ belongs to $\wo{E^c}$, then $\psi v$ belongs to
$\wo{{F\setminus E}}$. As $g=0$ in $\wm{{F\setminus E}}$
by~(\rf{pot}), it follows from~(\rf{psi}) that
$$
\int_{\rn} v\,d\lambda \Eq \langle g,\psi v\rangle \Eq 0 \,.
$$
In both cases $\int_{\rn}v\,d\lambda=0$. Since~$v$ and~$\lambda$ are
non-negative, this implies $v=0$ $\lambda$-\ae~$\rn$.

Similarly, it is enough to prove~(h) when~$v$ is non-negative. Let us fix
$v\in\wo{F}\cup\wo{F^c}$ with $v\ge 0$ \qe~$\rn$.
If~$v$ belongs to $\wo{F}$, then $(1-\psi) v$ belongs to
$\wo{{F\setminus E}}$. As $g=0$ in $\wm{{F\setminus E}}$
by~(\rf{pot}), it follows from~(\rf{psi}) that
$$
\int_{\rn} v\,d\nu \Eq -\langle g,(1-\psi) v\rangle \Eq 0 \,.
$$
If~$v$ belongs to $\wo{F^c}$, then $v=(1-\psi) v$ \qe~$\rn$.
Since ${g=0}$ in $\wm{F^c}$ by~(\rf{zero}), from~(\rf{psi}) we obtain
$$
\int_{\rn}v\,d\nu \Eq -\langle g, v\rangle \Eq 0 \,.
$$
In both cases $\int_{\rn}v\,d\nu=0$. Since~$v$ and~$\nu$ are
non-negative, this implies $v=0$ $\nu$-\ae~$\rn$.

It is enough to prove~(i) for every $\cp$-quasi open set~$U$ such that
either $U\subset E$ or $U\subset E^c$. In both cases by
Lemma~\rf{approx} there exists an increasing sequence $(v_j)$ of
functions of $\wo{E}\cup\wo{E^c}$, with
$0\le v_j\le 1^{}_U$ \qe~$\rn$, which converges to
$1^{}_U$ \qe~$\rn$. By~(g) we have $\int_{\rn} v_j\,d\lambda =0$ for
every~$j$, and passing to the limit as $j\to\infty$ we get
$\lambda(U)=0$.
The proof of~(j) is similar.

Since $(\partial E)^c\cap E$ and $(\partial E)^c\cap E^c$ are open sets,
by~(i) we have $\lambda\bigl((\partial E)^c\bigr)=0$, hence
$\supp\,\lambda\subset \partial E$. Similarly we prove that the
inclusion
$\supp\,\nu\subset \partial F$ follows from~(j).

Since $E$ is $\cp$-compatible with~$F$, the set ${E\setminus F}$ is
$\cp$-null. Consequently, by Remark~\rf{cap3}, the
$\ca$-potentials do not change if we replace~$E$ by ${E\cap F}$.
Therefore in the rest of the proof we may assume that $E\subset F$.

Let~$\chi$ be a function in $\cinftyo{\rn}$ such that ${\chi=1}$ in~$F$.
Since ${\chi-\psi}\in\wo{E^c}$, by~(f) and~(g) we have
${1-\psi}={\chi-\psi}=0$ $\lambda$-\ae~$\rn$. Since $\psi\in\wo{F}$,
by~(h) we have $\psi=0$ $\nu$-\ae~$\rn$. These facts imply that
$\lambda$ is concentrated in the set
$\{{\psi=1}\}$ while~$\nu$ is concentrated in the set $\{{\psi=0}\}$
and prove that~$\lambda$ and~$\nu$ are mutually singular.

Since $\psi=1$ $\lambda$-\ae~$\rn$ and $\psi=0$ \qe~$F^c$, by~(c) we
have
$\lambda(F^c)=0$. Similarly, as $\psi=0$
$\nu$-\ae~$\rn$ and $\psi=1$ \qe~$E$, by~(c) we have
$\nu(E)=0$.

Finally, condition~(b) determines uniquely the signed measure
$\lambda-\nu$. Since~$\lambda$ and~$\nu$ are non-negative and
mutually singular, by the uniqueness of the Hahn Decomposition we
have $\lambda=({\lambda-\nu})^+$ and $\nu=({\lambda-\nu})^-$, thus
the pair $(\lambda,\nu)$ is uniquely determined by conditions~(a)
and~(b). In particular~$\lambda$ and~$\nu$ do not depend on the
function~$\psi$ used in the proof.
\endproof
\defin{\thm{capdistr}}{Let~$E$ and~$F$ be two $\cp$-compatible
bounded sets in~$\rn$. The measures~$\lambda$ and~$\nu$ introduced
in the previous theorem are called the {\it inner\/} and the {\it outer
$\ca$-distributions\/} of~$E$ in~$F$.
}
\rem{\thm{wm}}{If $\overline E\subset \interior F$, it is easy to see
that the
$\ca$-distributions~$\lambda$ and~$\nu$ belong to $\wm{\rn}$. This is
not true, in general, when~$E$ and~$F$ are only $\cp$-compatible. For
a counterexample we refer to the Appendix of~[\rf{DM-GAR-cap}].
}
\prop{\thm{cd}}{Let~$E$ and~$F$ be two $\cp$-compatible
bounded sets in~$\rn$ and let~$\lambda$ and~$\nu$ be the inner and
the outer $\ca$-distributions of~$E$ in~$F$. Then
$$
\ca(E,F)= \lambda(\partial E)= \lambda(\rn)= \lambda(F)= \nu(\partial
F) =
\nu(\rn) = \nu(E^c) \,.
$$
}
\negbigskip
\proof By properties~(f) and~(k) of Theorem~\rf{capdistr1} we
have $\lambda(\partial E)= \lambda(\rn)= \lambda(F)$ and
$\nu(\partial F) =
\nu(\rn) = \nu(E^c)$.

Since $E$ is $\cp$-compatible with~$F$, the set ${E\setminus F}$ is
$\cp$-null. Consequently, by Remark~\rf{cap3}, the inner an the outer
$\ca$-distributions do not change if we replace~$E$ by ${E\cap F}$.
Therefore it is not restrictive to assume that $E\subset F$.
Let~$u$ be a $\ca$-potential of~$E$ in~$F$, let~$\psi$ be as in
Remark~\rf{compatible1}, and let~$\chi$ be a function in
$\cinftyo{\rn}$ such that $\chi=1$ in~$F$.
Since $\psi\in\wo{F}$ and ${\chi-\psi}\in\wo{E^c}$, by properties~(f),
(g), and (h) of Theorem~\rf{capdistr1} we have $\psi=0$
$\nu$-\ae~$\rn$, $\chi=1$ $\nu$-\ae~$\rn$, and
$\psi=\chi=1$ $\lambda$-\ae~$\rn$. By Theorem~\rf{capdistr1}(d) this
implies
$$
\displayliness{
\langle \A{u},\psi\rangle \Eq \int_{\rn} \psi \,d\lambda \Eq
\lambda(\rn) \,,
\cr
\langle \A{u},\chi-\psi\rangle \Eq -\int_{\rn} (\chi -\psi) \,d\nu
\Eq -\nu(\rn) \,.
\cr}\leqno(\frm{no})
$$
By Remark~\rf{cap1} we have
$\ca(E,F)=\langle \A{u},\psi\rangle$. Since
$Du=0$ \ae~$F^c$ and $D\chi=0$ \ae~$F$, by~(\rf{zero}) we have
$\langle \A{u},\chi-\psi\rangle = -\langle \A{u},\psi\rangle$.
Therefore~(\rf{no}) implies that
$\ca(E,F)=\lambda(\rn)=\nu(\rn)$.
\endproof
\parag{\chp{propcap}}{Monotonicity and continuity along monotone
sequences}
In this section we study the monotonicity and continuity properties of
$\ca(E,F)$ with respect to~$E$ and~$F$. These results are based
on the fundamental inequality proved in Lemma~\rf{fundamental}, on
the properties of the $\ca$-distributions discussed
in Section~\rf{capacity},  and on the properties of the minimal and
maximal $\ca$-potentials  proved in the following lemmas.
\lemma{\thm{comp1}}{Let~$E_1$, $E_2$,~$F$ be three bounded sets
in~$\rn$. Assume that $E_1\subset E_2$ and that~$E_2$ and~$F$ are
$\cp$-compatible. If~$u_1$ and~$u_2$ are the minimal
$\ca$-potentials of~$E_1$ and~$E_2$ in~$F$, then $u_1\le u_2$
\qe~$\rn$.
}
\proof By Corollary~\rf{restriction} $u_1$ coincides with the minimal
solution~$u$ of the problem
$$
\dirpb{u-u_1\in\wo{F\setminus E_2}}{\A{u}=0}{\wm{F\setminus E_2}}.
$$
By Remark~\rf{potential2} we have $u_1\le 1$ \qe~$E_2$. Therefore the
Comparison Principle (Theorem~\rf{comparison}) implies that $u_1\le
u_2$ \qe{} ${F\setminus E_2}$. Since $u_1\le 1=u_2$ \qe~$E_2$ and
$u_1=0=u_2$ \qe~$F^c$, we conclude that $u_1\le u_2$ \qe~$\rn$.
\endproof
\lemma{\thm{comp2}}{Let~$E$, $F_1$,~$F_2$ be three bounded sets
in~$\rn$. Assume that $F_1\subset F_2$ and that~$E$ and~$F_1$
are $\cp$-compatible. If~$u_1$ and~$u_2$ are the maximal
$\ca$-potentials of~$E$ in~$F_1$ and~$F_2$, then $u_1\le u_2$
\qe~$\rn$.
}
\proof By Corollary~\rf{restriction} $u_2$ coincides with the maximal
solution~$u$ of the problem
$$
\dirpb{u-u_2\in\wo{F_1\setminus E}}
{\A{u}=0}{\wm{F_1\setminus E}}.
$$
By Remark~\rf{potential2} we have $u_2\ge 0$ \qe~$F_1^c$. Therefore
the Comparison Principle (Theorem~\rf{comparison}) implies that
$u_1\le u_2$ \qe{} ${F_1\setminus E}$. Since $u_1= 1=u_2$ \qe~$E$ and
$u_1=0\le u_2$ \qe~$F_1^c$, we conclude that $u_1\le u_2$ \qe~$\rn$.
\endproof
\spazio
We prove now that the set function $\ca(\cdot,F)$ is increasing.
\th{\thm{monoE}}{Let~$E_1$, $E_2$,~$F$ be three bounded sets
in~$\rn$ such that $E_1\subset E_2$. Then $\ca(E_1,F)\le\ca(E_2,F)$.
}
\proof Since the inequality is trivial when~$E_2$ and~$F$ are not
$\cp$-compatible, the conclusion follows from Proposition~\rf{cd} and
from the following lemma.
\endproof
\lemma{\thm{monoE1}}{Let~$E_1$, $E_2$,~$F$ be three bounded sets
in~$\rn$. Assume that $E_1\subset E_2$ and that~$E_2$ and~$F$ are
$\cp$-compatible. Let~$\nu_1$
and~$\nu_2$ be the outer $\ca$-distributions of~$E_1$ and~$E_2$
in~$F$. Then $\nu_1(B)\le\nu_2(B)$ for every Borel set
${B\subset\rn}$.
}
\proof Let~$u_1$ and~$u_2$ be the minimal
$\ca$-potentials of~$E_1$ and~$E_2$ in~$F$. Then $u_1\le u_2$
\qe~$\rn$ by Lemma~\rf{comp1}. If we apply Lemma~\rf{fundamental}
with $B=F^c$, $C=E_2^c$, $w_1=u_1$, and $w_2=u_2$, we obtain
$\A{u_1}\ge\A{u_2}$ in $\wm{E_2^c}$. By~(d) and~(g) of
Theorem~\rf{capdistr1}
we have
$$
\int_{\rn} v\,d\nu_1 \Eq -\langle \A{u_1}, v\rangle \Le
-\langle \A{u_2}, v\rangle \Eq \int_{\rn} v\,d\nu_2
$$
for every $v\in\wo{E_2^c}$ with $v\ge0$ \qe~$\rn$. By
Lemma~\rf{approx} this implies $\nu_1(V)\le\nu_2(V)$ for
every $\cp$-quasi open set~$V$ contained in~$E_2^c$. As $\{{u_2>0}\}$
is $\cp$-quasi open, by Theorem~\rf{capdistr1}(j) we have
$\nu_1(\{{u_2\ge1}\})\le
\nu_1(\{{u_2> 0}\})=0$ and $\nu_2(\{{u_2\ge1}\})\le
\nu_2(\{{u_2> 0}\})=0$. For every open set $U\subset\rn$
the set $U\cap\{{u_2<1}\}$ is $\cp$-quasi open and is contained
in~$E_2^c$ (up to a $\cp$-null set). Therefore
$\nu_1(U)=\nu_1(U\cap\{{u_2<1}\})\le
\nu_2(U\cap\{{u_2<1}\})=\nu_2(U)$. Since~$\nu_1$
and~$\nu_2$ are Radon measures, this implies that
$\nu_1(B)\le\nu_2(B)$ for every Borel set ${B\subset\rn}$.
\endproof
\spazio
The following theorem shows that $\ca(E,\cdot)$ is decreasing.
\th{\thm{monoF}}{Let~$E$, $F_1$,~$F_2$ be three bounded sets
in~$\rn$ such that $F_1\subset F_2$. Then $\ca(E, F_1)\ge\ca(E,F_2)$.
}
\proof Since the inequality is trivial when~$E$ and~$F_1$ are not
$\cp$-compatible, the conclusion follows from Proposition~\rf{cd} and
from the
following lemma.
\endproof
\lemma{\thm{monoF1}}{Let~$E$, $F_1$,~$F_2$ be three bounded sets
in~$\rn$. Assume that $F_1\subset F_2$ and that~$E$ and~$F_1$ are
$\cp$-compatible. Let~$\lambda_1$
and~$\lambda_2$ be the inner $\ca$-distributions of~$E$ in~$F_1$
and~$F_2$. Then $\lambda_1(B)\ge\lambda_2(B)$ for every Borel set
${B\subset\rn}$.
}
\proof Let~$u_1$ and~$u_2$ be the maximal
$\ca$-potentials of~$E$ in~$F_1$ and~$F_2$. Then $u_1\le u_2$
\qe~$\rn$ by Lemma~\rf{comp2}. If we apply Lemma~\rf{fundamental}
with $B=E$, $C=F_1$, $w_1=u_1$, and $w_2=u_2$, we obtain
$\A{u_1}\ge\A{u_2}$ in $\wm{F_1}$. By~(d) and~(h) of
Theorem~\rf{capdistr1}
we have
$$
\int_{\rn} v\,d\lambda_1 \Eq \langle \A{u_1}, v\rangle \Ge
\langle \A{u_2}, v\rangle \Eq \int_{\rn} v\,d\lambda_2
$$
for every $v\in\wo{F_1}$ with $v\ge0$ \qe~$\rn$. By
Lemma~\rf{approx} this implies $\lambda_1(V)\ge\lambda_2(V)$ for
every $\cp$-quasi open set~$V$ contained in~$F_1$. As $\{{u_1<1}\}$ is
$\cp$-quasi open, by  Theorem~\rf{capdistr1}(i) we have
$\lambda_1(\{{u_1\le0}\})\le
\lambda_1(\{{u_1<1}\})=0$ and  $\lambda_2(\{{u_1\le0}\})\le
\lambda_2(\{{u_1<1}\})=0$. For every open set $U\subset\rn$ the set
$U\cap\{{u_1>0}\}$ is $\cp$-quasi open and is contained in~$F_1$ (up to
a $\cp$-null set).
Therefore $\lambda_1(U)=\lambda_1(U\cap\{{u_1>0}\})\ge
\lambda_2(U\cap\{{u_1>0}\})=\lambda_2(U)$. Since~$\lambda_1$
and~$\lambda_2$ are Radon measures, this implies that
$\lambda_1(B)\ge\lambda_2(B)$ for every Borel set ${B\subset\rn}$.
\endproof
\spazio
The following theorem shows that the set function $\ca(\cdot,F)$ is
continuous along all increasing sequences.
\th{\thm{increasingE}}{Let~$E$ and~$F$ be two bounded sets in~$\rn$.
If~$E$ is the union of an increasing sequence of sets $(E_j)$, then
$$
\ca(E,F)\Eq \lim_{j\to\infty} \ca(E_j,F) \Eq \sup_{j} \ca(E_j,F)\,.
$$
}
\negbigskip
\proof Let $S=\sup_{j} \ca(E_j,F)$. By monotonicity
(Theorem~\rf{monoE}) we have $S\le \ca(E,F)$. It remains to prove the
opposite inequality when ${S\finite}$, and hence each
set~$E_j$ is $\cp$-compatible with~$F$. Let~$u$ and~$u_j$ be
the minimal $\ca$-potentials of~$E$ and~$E_j$ in~$F$.
As ${S\finite}$, by~(\rf{coerc1}) the sequence $(u_j)$ is bounded in
$\wp{\rn}$, and by~(\rf{bound1}) the sequence $(\A{u_j})$ is
bounded in  $\wm{\rn}$. Passing, if necessary, to a
subsequence, we may assume  that $(u_j)$ converges weakly
in $\wp{\rn}$ to some function $w\in\wo{F}$ and  that  $(\A{u_j})$
converges weakly in $\wm{\rn}$ to some element~$f$ of $\wm{\rn}$.

We want to prove that $\A{u}=f$ in $\wm{F}$ and that $w=u$
\qe~$\rn$.
 From the monotonicity condition~(\rf{mono1}) for every~$j$ we
obtain
$$
\langle \A{v}, v-u_j \rangle \Ge \langle \A{u_j}, v-u_j \rangle
\qquad \forall v\in\wp{\rn}\,.
\leqno(\frm{200})
$$
If $j\le i$, by Lemma~\rf{comp1} we have $u_j\le u_i\le u$ \qe~$\rn$,
hence $u_j\le w\le u$ \qe~$\rn$ for every~$j$.
Since $\A{u_j}\ge0$ in $\wm{F}$
(Lemma~\rf{signA}), we have
$$
\langle \A{u_j}, v-u_j \rangle
\Ge \langle \A{u_j}, v-w \rangle \qquad\forall v\in\wo{F}\,,
$$
which, together with~(\rf{200}), gives
$$
\langle \A{v}, v-u_j \rangle \Ge \langle \A{u_j}, v-w\rangle
\qquad \forall v\in\wo{F} \,.
$$
Passing to the limit as $j\to\infty$ we get
$$
\langle \A{v}, v-w\rangle \Ge \langle f , v-w\rangle
\qquad\forall v\in \wo{F}\,.
$$
Putting $v=w+\e z$, with $z\in\wo{F}$ and ${\e>0}$, and dividing
by~$\e$ we obtain
$$
\langle A(w+\e z), z \rangle \Ge \langle f , z\rangle
\qquad\forall z\in \wo{F}\,.
$$
Passing to the limit as $\e\to0$ we get
$$
\langle \A{w}, z\rangle \Ge \langle f , z\rangle
\qquad\forall z \in \wo{F}\,,
$$
hence $\A{w}=f$ in $\wm{F}$. As $u_j\le w\le u$ \qe~$\rn$ and
$u_j=u=1$ \qe~$E_j$ (Remark~\rf{potential1}), we have $w=1$
\qe~$E_j$ for every~$j$, hence $w=1$
\qe~$E$. Since $w\in\wo{F}$, we have also $w=0$ \qe~$F^c$. This shows
that~$w$ satisfies the first condition in~(\rf{pot}).

It remains to prove that $\A{w}=0$ in $\wm{{F\setminus E}}$. If
$v\in\wo{{F\setminus E}}$, then $v\in \wo{{F\setminus E_j}}$
and the definition of~$u_j$ implies that
$\langle \A{u_j}, v \rangle \Eq0$
for every~$j$. Since $(\A{u_j})$ converges to~$f$ weakly in
$\wm{\rn}$ and $\A{w}=f$ in $\wm{F}$, we conclude that
$$
\langle \A{w},v\rangle \Eq \langle f,v\rangle
\Eq \lim_{j\to\infty} \langle \A{u_j},v\rangle  \Eq0
\qquad \forall v\in\wo{{F\setminus E}}\,,
$$
hence $\A{w}=0$ in $\wm{{F\setminus E}}$.
This proves that $w$ is a $\ca$-potential of~$E$ in~$F$. Since ${w\le u}$
\qe~$\rn$, by the minimality of~$u$ we obtain ${w=u}$ \qe~$\rn$.

By Remark~\rf{cap1} we have $\ca(E_j,F)=\langle \A{u_j},u\rangle$ for
every~$j$.
Using again the fact that $(\A{u_j})$ converges to~$f$ weakly in
$\wm{\rn}$ and that $\A{u}=f$ in $\wm{F}$ we obtain
$$
\ca(E,F)\Eq \langle \A{u},u\rangle  \Eq \langle f,u\rangle
\Eq \lim_{j\to\infty} \langle \A{u_j},u\rangle
\Eq \lim_{j\to\infty} \ca(E_j,F)\,,
$$
which concludes the proof of the theorem.
\endproof
\spazio
For the continuity of  the set function $\ca(\cdot,F)$ along decreasing
sequences $(E_j)$ we need two additional assumptions: the sets $E_j$
must be $\cp$-quasi closed and $\cp$-compatible with~$F$.
\th{\thm{decreasingE}}{Let~$F$ be a bounded set in~$\rn$, let $(E_j)$ be
a decreasing sequence of $\cp$-quasi closed bounded sets, and let~$E$
be their intersection. If~$E_1$ and~$F$ are $\cp$-compatible, then
$$
\ca(E,F)\Eq \lim_{j\to\infty} \ca(E_j,F) \Eq \inf_{j} \ca(E_j,F)\,.
$$
}
\negbigskip
\proof By Theorem~\rf{monoE} we have $\ca(E_j,F)\le\ca(E_1,F)\finite$
for every~$j$. Let~$u_j$ be the minimal $\ca$-potential of~$E_j$
in~$F$. By~(\rf{coerc1}) the sequence $(u_j)$ is bounded in
$\wp{\rn}$, and by~(\rf{bound1}) the sequence $(\A{u_j})$ is
bounded in  $\wm{\rn}$. Passing, if necessary, to a
subsequence, we may assume  that $(u_j)$ converges weakly
in $\wp{\rn}$ to some function $u\in\wo{F}$ and  that $(\A{u_j})$
converges weakly in $\wm{\rn}$ to some element~$f$ of $\wm{\rn}$.

We want to prove that $\A{u}=f$ in $\wm{F}$ and that~$u$ is a
$\ca$-potential of~$E$ in~$F$.
 From the monotonicity condition~(\rf{mono1}) for every~$j$ we
obtain
$$
\langle \A{v}, v-u_j \rangle \Ge \langle \A{u_j}, v-u_j \rangle
\qquad \forall v\in\wp{\rn}\,.
\leqno(\frm{202})
$$
If $j\ge i$, by Lemma~\rf{comp1} we have $u_j\le u_i$ \qe~$\rn$.
Since $\A{u_j}\ge0$ in $\wm{F}$
(Lemma~\rf{signA}), we have
$$
\langle \A{u_j}, v-u_j \rangle
\Ge \langle \A{u_j}, v-u_i \rangle \qquad\forall v\in\wo{F}\,,
$$
which, together with~(\rf{202}), gives
$$
\langle \A{v}, v-u_j \rangle \Ge \langle \A{u_j}, v-u_i\rangle
\qquad \forall v\in\wo{F}
$$
whenever $j\ge i$. Passing to the limit as $j\to\infty$ we obtain
$$
\langle \A{v}, v-u\rangle \Ge \langle f , v-u_i\rangle
\qquad\forall v\in \wo{F}\,,
$$
and as $i\to\infty$ we get
$$
\langle \A{v}, v-u\rangle \Ge \langle f , v-u\rangle
\qquad\forall v\in \wo{F}\,.
$$
Putting $v=u+\e z$, with $z\in\wo{F}$ and ${\e>0}$, and dividing
by~$\e$ we obtain
$$
\langle A(u+\e z), z \rangle \Ge \langle f , z\rangle
\qquad\forall z\in \wo{F}\,.
$$
Passing to the limit as $\e\to0$ we get
$$
\langle \A{u}, z\rangle \Ge \langle f , z\rangle
\qquad\forall z \in \wo{F}\,,
$$
hence $\A{u}=f$ in $\wm{F}$.
Since $u_j=1$ \qe~$E$ and $u_j=0$ \qe~$F^c$ for every~$j$
(Remark~\rf{potential1}), we have
$u=1$ \qe~$E$ and $u=0$ \qe~$F^c$. This shows that~$u$ satisfies the
first condition in~(\rf{pot}).

It remains to prove that $\A{u}=0$ in $\wm{{F\setminus E}}$.
Let us fix $v\in \wm{{F\setminus E}}$. Since the sets $E_j$ are
$\cp$-quasi closed, by Lemma~\rf{approx2} there exists a sequence
$(v_j)$ which converges to~$v$ strongly in $\wp{\rn}$ and such that
$v_j\in\wo{F\setminus E_j}$ for every~$j$. As the sequence $(E_j)$ is
decreasing, we have $v_i\in\wo{F\setminus E_j}$ for every ${j\ge i}$.
By the definition of~$u_j$ we have $\langle \A{u_j},v_i\rangle=0$ for
every ${j\ge i}$. Since $(\A{u_j})$ converges to~$f$ in $\wm{\rn}$ and
$\A{u}=f$ in $\wm{F\setminus E}$, as ${j\to\infty}$ we get
$$
\langle \A{u},v_i\rangle \Eq \langle f,v_i\rangle \Eq
\lim_{j\to\infty} \langle \A{u_j},v_i\rangle \Eq0 \,.
$$
Passing to the limit as ${i\to\infty}$ we obtain
$\langle \A{u},v\rangle=0$, hence $\A{u}=0$ in $\wm{F\setminus E}$.
This proves that $u$ is a $\ca$-potential of~$E$ in~$F$.

By Remark~\rf{cap1} we have $\ca(E,F)=\langle \A{u},u_1\rangle$ and
$\ca(E_j,F)=\langle \A{u_j},u_1\rangle$ for every~$j$.
Using again the fact that $(\A{u_j})$ converges to~$f$ weakly in
$\wm{\rn}$ and that $\A{u}=f$ in $\wm{F}$ we obtain
$$
\ca(E,F)\Eq \langle \A{u},u_1\rangle  \Eq \langle f,u_1\rangle
\Eq \lim_{j\to\infty} \langle \A{u_j},u_1\rangle
\Eq \lim_{j\to\infty} \ca(E_j,F)\,,
$$
which concludes the proof of the theorem.
\endproof
\rem{\thm{decreasing1}}{Elementary examples in the case ${p=2}$ and
$\A{u}=-\Delta u$ show that the conclusion of Theorem~\rf{decreasingE}
does not hold if the sets~$E_j$ are not $\cp$-quasi closed. The
assumption that~$E_1$ and~$F$ are $\cp$-compatible is automatically
satisfied if~$F$ is open and the sets~$E_j$ are compact and contained
in~$F$.
}
\spazio
We consider now the continuity properties with respect to~$F$. The
following theorem shows that the set function $\ca(E,\cdot)$ is
continuous along all decreasing sequences.
\th{\thm{decreasingF}}{Let~$E$ and~$F$ be two bounded sets in~$\rn$.
If~$F$ is the intersection of a decreasing sequence of sets $(F_j)$, then
$$
\ca(E,F)\Eq \lim_{j\to\infty} \ca(E,F_j) \Eq \sup_{j} \ca(E,F_j)\,.
$$
}
\negbigskip
\proof It is enough to repeat the proof of Theorem~\rf{increasingE}
with obvious modifications. For instance we have
to replace the minimal
$\ca$-potentials by the maximal $\ca$-potentials, $\wo{F}$ by
$\{{u\in\wp{\rn}}:{u=1\,\ \cp\hbox{-q.e.\ in\ }\,E}\}$, and $\wm{F}$ by
$\wm{E^c}$.
\endproof
\spazio
For the continuity of  the set function $\ca(E,\cdot)$ along increasing
sequences, we need additional assumptions.
\th{\thm{increasingF}}{Let~$E$ and~$F$ be two bounded set in~$\rn$.
Assume that $F$ is the union of an increasing sequence $(F_j)$ of
$\cp$-quasi open sets such that~$E$ and~$F_1$ are $\cp$-compatible.
Then
$$
\ca(E,F)\Eq \lim_{j\to\infty} \ca(E,F_j) \Eq \inf_{j} \ca(E,F_j)\,.
$$
}
\negbigskip
\proof It is enough to modify the proof of Theorem~\rf{decreasingE} as
in
the proof of Theorem~\rf{decreasingF}.
\endproof
\parag{\chp{app}}{Approximation properties and subadditivity}
In this section we prove that, if~$E$ and~$F$ are bounded Borel sets,
then $\ca(E,F)$ can be approximated by $\ca(K,U)$, with~$K$ compact,
$K\subset E$, and~$U$ bounded and open, $U\supset F$. Finally we
prove that
$\ca(E,F)$ is countably subadditive with respect to~$E$.

We begin with the problem of the approximation of $\ca(E,F)$ by
$\ca(K,F)$, with~$K$ compact, $K\subset E$.
\lemma{\thm{choqE}}{Let~$E$ and~$F$ be two
$\cp$-compatible bounded sets in~$\rn$. If~$E$ is a Borel set, then
$$
\ca(E,F)\Eq
\sup\{\ca(K,F): K\,\ \hbox{\rm compact,\ }\, K\subset E\} \,.
\leqno(\frm{c1})
$$
}
\negbigskip
\proof Let~$\psi$ be a function in $\wo{F}$ such that $\psi=1$ \qe~$E$,
let $H=\{{\psi\ge1}\}$, and let $\alpha$ be the set function defined by
$\alpha(B)=\ca(B\cap H,F)$ for every $B\subset\rn$. Since~$H$ is
$\cp$-quasi closed and $\cp$-compatible with~$F$, the set
function~$\alpha$ satisfies the following properties:
\smallskip
\item{(i)} if $B\subset C$, then $\alpha(B)\le\alpha(C)$
(Theorem~\rf{monoE});
\smallskip
\item{(ii)} if~$B$ is the union of an increasing sequence of sets $(B_j)$,
then $\alpha(B)=\sup_j \alpha(B_j)$
(Theorem~\rf{increasingE});
\smallskip
\item{(iii)} if~$K$ is the intersection of a decreasing sequence of compact
sets $(K_j)$, then $\alpha(K)=\inf_j \alpha(K_j)$
(Theorem~\rf{decreasingE}).
\par\noindent
Therefore~$\alpha$ is an abstract capacity in the sense of Choquet. By
the Capacitability Theorem ([\rf{CHO}], Theorem~1) for every Borel set
$B\subset\rn$ we have
$$
\alpha(B) \Eq
\sup\{\alpha(K): K\,\ \hbox{compact,\ }\, K\subset B\} \,.
\leqno(\frm{c2})
$$
Since $\psi=1$ \qe~$E$, we have $\alpha(B)=\ca(B,F)$ for every
$B\subset E$ (Remark~\rf{cap3}). Consequently~(\rf{c2})
implies~(\rf{c1}).
\endproof
\th{\thm{choquetE}}{Let~$E$ and~$F$ be two bounded sets in~$\rn$.
If~$E$ is a Borel set, then
$$
\ca(E,F)\Eq
\sup\{\ca(K,F): K\,\ \hbox{\rm compact,\ } \, K\subset E \} \,.
\leqno(\frm{c3})
$$
}
\negbigskip
\proof Let~$\Om$ be a bounded open set containing~$E$ and~$F$,
let~$D$ be a
countable dense subset of $\wo{F}$, and let~$F_0$ be the union of the
sets $\{{v\neq 0}\}$ for ${v\in D}$. If $u\in\wo{F}$, then there exists a
sequence $(v_j)$ in~$D$ which converges to~$u$ strongly in $\wp{\rn}$.
Since $v_j=0$ \qe~$F_0^c$, we have $u=0$ \qe~$F_0^c$. Therefore,
if~$B$ is a bounded set in~$\rn$ such that $\cp({B\setminus
F_0},\Om)>0$, then~$B$ and~$F$ are not $\cp$-compatible. We may
assume that all functions $v\in D$ are Borel functions, so that~$F_0$ is a
Borel set.

If $\cp({E\setminus F_0},\Om)>0$, then $\ca(E,F)=\pinfty$ by the
previous remark. If we apply
Choquet's Capacitability Theorem ([\rf{CHO}], Theorem~1) to the capacity
$\cp(\cdot,\Om)$, we obtain that there exists a compact set~$K$
contained in ${E\setminus F_0}$ such that
$\cp({K\setminus F_0},\Om)=\cp(K,\Om)>0$. This implies that
$\ca(K,F)=\pinfty$ and proves~(\rf{c3}) in this case.

If $\cp({E\setminus F_0},\Om)=0$, then $\ca(E,F)=\ca({E\cap F_0},F)$
(Remark~\rf{cap3}) and ${E\cap F_0}$ is the union of the sets
$E\cap\{{|v|>{1\over k}}\}$ for $k\in\n$ and $v\in D$. Since all these
sets are
$\cp$-compatible with~$F$ (with $\psi=(k|v|)\land1$), so are their finite
unions
(Remark~\rf{compatible3}). Therefore ${E\cap F_0}$ is the union of an
increasing sequence $(E_j)$ of sets which are $\cp$-compatible
with~$F$. By
Theorem~\rf{increasingE} we have
$$
\ca(E,F)\Eq\ca({E\cap F_0},F)\Eq \sup_j \ca(E_j,F)\,,
$$
and~(\rf{c3}) follows from the fact that
$$
\ca(E_j,F)\Eq
\sup\{\ca(K,F): K\,\ \hbox{compact,\ } \, K\subset E_j \}
$$
by Lemma~\rf{choqE}.
\endproof
\spazio
We consider now the problem of the approximation of $\ca(E,F)$ by
$\ca(E,U)$, with~$U$ bounded and open, ${U\supset F}$.
\lemma{\thm{choqF}}{Let~$E$ and~$F$ be two
$\cp$-compatible bounded sets in~$\rn$. If~$F$ is a Borel set, then
$$
\ca(E,F)\Eq
\sup\{\ca(E,U): U\,\ \hbox{\rm bounded\ and\ open,\ }\,
U\supset F\} \,.
\leqno(\frm{c11})
$$
}
\negbigskip
\proof Let~$\Om$ be a bounded open set in~$\rn$ containing~$F$,
let~$\psi$ be a function in $\wo{F}$ such that $\psi=1$ \qe~$E$, let
$V=\{{\psi>0}\}$, and let $\beta$ be the set function defined by
$\beta(B)=\ca(E,{(\Om\setminus B)\cup V})$ for every $B\subset\rn$.
Since~$V$ is $\cp$-quasi open and~$E$ is $\cp$-compatible with~$V$,
the set function~$\beta$ satisfies the following properties:
\smallskip
\item{(i)} if $B\subset C$, then $\beta(B)\le\beta(C)$
(Theorem~\rf{monoF});
\smallskip
\item{(ii)} if~$B$ is the union of an increasing sequence of sets $(B_j)$,
then $\beta(B)=\sup_j \beta(B_j)$
(Theorem~\rf{decreasingF});
\smallskip
\item{(iii)} if~$K$ is the intersection of a decreasing sequence of compact
sets $(K_j)$, then $\beta(K)=\inf_j \beta(K_j)$
(Theorem~\rf{increasingF}).
\par\noindent
Therefore~$\beta$ is an abstract capacity in the sense of Choquet. By the
Capacitability Theorem ([\rf{CHO}], Theorem~1) for every Borel set
$B\subset\rn$ we have
$$
\beta(B) \Eq
\sup\{\beta(K): K\,\ \hbox{compact,\ }\, K\subset B\} \,.
\leqno(\frm{c12})
$$
Since $\psi=0$ \qe~$F^c$, we have $\beta(B)=\ca(E,{\Om\setminus B})$
for every $B\subset F^c$ (Remark~\rf{cap3}). In particular
$\beta(F^c)=\ca(E,F)$.
Consequently~(\rf{c12}) gives
$$
\ca(E,F)\Eq
\sup\{\ca(E,\Om\setminus K): K\,\ \hbox{compact,\ }\, K\subset F^c\}
\,,
$$
which implies~(\rf{c11}).
\endproof
\th{\thm{choquetF}}{Let~$E$ and~$F$ be two bounded sets in~$\rn$.
If~$F$ is a Borel set, then
$$
\ca(E,F)\Eq
\sup\{\ca(E,U): U\,\ \hbox{\rm bounded\ and\ open,\ } \, U\supset F
\} \,.
\leqno(\frm{c13})
$$
}
\negbigskip
\proof Let~$\Om$ be a bounded open set containing~$\overline E$
and~$\overline F$,
let~$D$ be a countable dense subset of
$H=\{{v\in\wp{\rn}}:{v\ge1\,\ \cp\hbox{-q.e.\ in\ }\,E}\}$, and
let~$E_0$ be the intersection of the sets $\{{v\ge1}\}$ for ${v\in D}$. If
$u\in H$, then there exists a sequence $(v_j)$ in~$D$ which converges
to~$u$ strongly in $\wp{\rn}$.
Since $v_j\ge1$ \qe~$E_0$, we have $u\ge1$ \qe~$E_0$. Therefore,
if~$B$
is a bounded set in~$\rn$ such that $\cp({E_0\setminus B},\Om)>0$,
then~$E$ and~$B$ are not $\cp$-compatible. We may
assume that all functions $v\in D$ are Borel functions, so that~$E_0$ is a
Borel set. As $\overline E\subset\Om$, the set ${E_0\setminus \Om}$ is
$\cp$-null, thus we may assume that ${E_0\subset\Om}$.

If $\cp({E_0\setminus F},\Om)>0$, then $\ca(E,F)=\pinfty$ by the
previous remark. If we apply Choquet's Capacitability Theorem
([\rf{CHO}], Theorem~1) to the capacity $\cp(\cdot,\Om)$, we obtain that
there exists a compact set~$K$ contained in ${E_0\setminus F}$ such
that $\cp(K,\Om)>0$. As
$\cp({E_0\setminus (\Om\setminus K)}, \Om)= \cp(K,\Om)>0$, we obtain
that
$\ca(E,{\Om\setminus K})=\pinfty$ and~(\rf{c13}) is proved.

If $\cp({E_0\setminus F},\Om)=0$, then $\ca(E,F)=\ca(E,{F\cup E_0})$
(Remark~\rf{cap3}) and ${F\cup E_0}$ is the intersection of the sets
$F\cup\{{v>1-{1\over k}}\}$ for $k\in\n$ and $v\in D$. Since~$E$ is
$\cp$-compatible with all these sets (with $\psi=(kv-k+1)^+$),
$E$ is $\cp$-compatible with their finite
intersections
(Remark~\rf{compatible3}). Therefore ${F\cup E_0}$ is the intersection
of a decreasing sequence $(F_j)$ of sets  such that~$E$ is
$\cp$-compatible with~$F_j$ for every~$j$. By
Theorem~\rf{decreasingF} we have
$$
\ca(E,F)\Eq\ca({E},F\cup E_0)\Eq \sup_j \ca(E,F_j)\,,
$$
and~(\rf{c13}) follows from the fact that
$$
\ca(E,F_j)\Eq
\sup\{\ca(E,U): U\,\ \hbox{bounded\ and\ open,\ } \, U\supset F_j \}
$$
by Lemma~\rf{choqF}.
\endproof
\spazio
We are now in a position to prove the main approximation theorem
for~$\ca$.
\th{\thm{choquet}}{Let~$E$ and~$F$ be two bounded Borel sets
in~$\rn$. Then
$$
\ca(E,F)\Eq
\sup\{\ca(K,U): K\,\ \hbox{\rm compact,\ } \, K\subset E\,,
\,\ U\,\ \hbox{\rm bounded\ and\ open,\ } \, U\supset F \} \,.
$$
}
\negbigskip
\proof The conclusion follows from Theorems~\rf{choquetE}
and~\rf{choquetF}.
\endproof
\spazio
We consider now the problem of the approximation of $\ca(E,F)$ by
$\ca(U,F)$, with~$U$ bounded and $\cp$-quasi open, ${U\supset E}$.
\prop{\thm{qopen}}{Let~$E$ and~$F$ be two bounded sets in~$\rn$.
Then
$$
\ca(E,F) \Eq
\inf\{\ca(U,F): U\,\ \hbox{\rm bounded\ and\ }\,\cp\hbox{\rm-quasi\
open,\ } \,
U\supset E\} \,.
\leqno(\frm{aqo})
$$
}
\negbigskip
\proof Let~$I$ be the right hand side of~(\rf{aqo}). By monotonicity we
have $\ca(E,F)\le I$. Let us prove the opposite inequality when
$\ca(E,F)\finite$, and hence~$E$ and~$F$ are $\cp$-compatible. Let~$u$
be the minimal $\ca$-potential of~$E$ in~$F$. For every
$k\in\n$ let $E_k=\{{u\ge 1-{1\over 2k}}\}$,
$U_k= \{{u> 1-{1\over 2k}}\}$, and let $E_0=\{{u=1}\}$.  By
Remark~\rf{potential2} we can write $E_0=\{{u\ge1}\}$, hence $E_0$ is
the intersection of the decreasing sequence $(E_k)$. It is easy to
see that~$u$ is a
$\ca$-potential of~$E_0$ in~$F$, hence
$$
\ca(E_0,F)\Eq \iaa{F}{Du}{Du} \Eq \ca(E,F)\,.
$$
Since~$E_1$ and~$F$ are
$\cp$-compatible, with $\psi={(2u)\land1}$, and all sets~$E_k$ are
$\cp$-quasi closed,  by Theorem~\rf{decreasingE} we
have
$$
\ca(E,F)\Eq \ca(E_0,F)\Eq \inf_{k} \ca(E_k,F)\,.
$$
Since the sets $U_k$ are $\cp$-quasi open and
$E\subset U_k\subset E_k$ (up to a $\cp$-null set), we have
$$
\ca(E,F) \Eq \inf_{k} \ca(E_k,F) \Ge \inf_{k} \ca(U_k,F) \Ge I\,,
$$
which concludes the proof.
\endproof
\spazio
We are now in a position to prove the subadditivity of the
capacity~$\ca(\cdot,F)$.
\th{\thm{subadd}}{Let~$E_1$, $E_2$,~$F$ be three bounded sets
in~$\rn$. Then
$$
\ca(E_1\cup E_2,F) \le
\ca(E_1,F) + \ca(E_2,F) \,.
\leqno(\frm{s0})
$$
}
\negbigskip
\proof The inequality is trivial if $\ca(E_1,F)=\pinfty$ or
$\ca(E_2,F)=\pinfty$. Therefore we may assume that~$E_1$ and~$E_2$
are $\cp$-compatible with~$F$. By Proposition~\rf{qopen} for every
${\e>0}$ there exist two $\cp$-quasi open sets~$U_1$ and~$U_2$ such
that $E_1\subset U_1$, $E_2\subset U_2$, and
$$
\ca(U_1,F) + \ca(U_2,F) \le \ca(E_1,F) + \ca(E_2,F) +\e\,.
\leqno(\frm{s1})
$$
By Lemma~\rf{compact} there exist two increasing sequences
$(K_1^j)$ and $(K_2^j)$ of compact sets, contained in~$U_1$ and
${U_2\setminus U_1}$ respectively, whose unions cover $\cp$-quasi all
of~$U_1$ and ${U_2\setminus U_1}$. By
Remark~\rf{cap3} and by Theorems~\rf{monoE} and~\rf{increasingE} we
have
$$
\ca(E_1\cup E_2,F)\Le \ca(U_1\cup U_2,F) \Eq
\lim_{j\to\infty} \ca(K_1^j\cup K_2^j,F) \,.
\leqno(\frm{s2})
$$
Since by monotonicity (Theorem~\rf{monoE})
$$
\ca(K_1^j,F) + \ca(K_2^j,F) \le \ca(U_1,F) + \ca(U_2,F)\,,
$$
in view of~(\rf{s1}) and~(\rf{s2}) it is enough to prove that for every~$j$
we have
$$
\ca(K_1^j\cup K_2^j,F)\le \ca(K_1^j,F) + \ca(K_2^j,F)\,.
\leqno(\frm{s3})
$$
Let us fix~$j$ and let $K_1=K_1^j$ and $K_2=K_2^j$. As the
compact sets $K_1$ and $K_2$ are disjoint, there exist
two disjoint open set $V_1$ and $V_2$ such that
$K_1\subset V_1$ and $K_2\subset V_2$.

Let~$u_1$, $u_2$, and~$u$ be the minimal $\ca$-potentials of~$K_1$,
$K_2$, and ${K_1\cup K_2}$ in~$F$, and let~$\lambda_1$, $\lambda_2$,
and~$\lambda$ be the corresponding inner $\ca$-distributions. We
want to
prove that
$$
\lambda(B\cap K_1)\le\lambda_1(B)\quad \hbox{and} \quad
\lambda(B\cap K_2)\le\lambda_2(B)
\quad \hbox{for\ every\ Borel\ set\ }\,  B\subset \rn\,.
\leqno(\frm{s4})
$$
By Lemma~\rf{comp1} we have $u_1\le u$ and $u_2\le u$ \qe~$\rn$.
If we apply Lemma~\rf{fundamental} with $B=K_1$, $C=V_1\cap F$,
$w_1=u_1$, $w_2=u$, we obtain $\A{u}\le\A{u_1}$ in
$\wm{{V_1\cap F}}$. By properties~(d) and~(h) of
Theorem~\rf{capdistr1} we have
$$
\int_{\rn} v\,d\lambda \Eq \langle \A{u}, v\rangle \Le
\langle \A{u_1}, v\rangle  \Eq \int_{\rn} v\,d\lambda_1
$$
for every $v\in\wo{{V_1\cap F}}$ with $v\ge0$ \qe~$\rn$. By
Lemma~\rf{approx} this implies $\lambda(V)\le\lambda_1(V)$ for
every
$\cp$-quasi open set~$V$ in ${V_1\cap F}$. In particular, since
$\{{u_1>0}\}$ is $\cp$-quasi open and is contained in~$F$, we have
$\lambda(U\cap V_1 \cap\{{u_1{>}0}\})\le \lambda_1(U\cap V_1 \cap
\{{u_1{>}0}\})$ for every open
set~$U$ in~$\rn$. As~$\lambda$ and~$\lambda_1$ are Radon measures,
this
implies that $\lambda(B)\le \lambda_1(B)$ for every Borel set
$B\subset {V_1 \cap \{{u_1>0}\}}$. Since $u_1=1$ \qe~$K_1$ and
$K_1\subset V_1$, we have
$\lambda(B\cap K_1)\le \lambda_1(B\cap K_1)\le \lambda_1(B)$ for
every Borel set
$B\subset\rn$, which proves the first inequality in~(\rf{s4}). The
other inequality is proved in a similar way.

Since $\supp\,\lambda\subset K_1\cup K_2$
(Theorem~\rf{capdistr1}(f)),
by~(\rf{s4}) we have $\lambda(\rn)=\lambda(K_1)+\lambda(K_2)\le
\lambda_1(\rn)+\lambda_2(\rn)$, which gives~(\rf{s3}) by
Proposition~\rf{cd}.
\endproof
\spazio
When~$\Om$ is a bounded open set in~$\rn$ and $E\subset \Om$,
then $\ca(E,\Om)$ can be approximated by $\ca(U,\Om)$, with~$U$
open, $E\subset U\subset\Om$.
\prop{\thm{infopen}}{Let~$\Om$ be a bounded open set in~$\rn$ and
let $E\subset \Om$. Then
$$
\ca(E,\Om) \Eq
\inf\{\ca(U,\Om): U\,\ \hbox{\rm open,\ }\, E\subset U\subset
\Om\} \,.
\leqno(\frm{infop})
$$
}
\negbigskip
\proof Let~$I$ be the right hand side of~(\rf{infop}). Since
$\ca(E,F)\le I$ by monotonicity,  we have only to prove the opposite
inequality. By Theorem~\rf{monoE} and by Proposition~\rf{qopen} for
every ${\e>0}$ there exists a $\cp$-quasi open set~$V$ such that
$E\subset
V\subset\Om$ and
$$
\ca(V,\Om) \Le \ca(E,\Om)+\e \,.
\leqno(\frm{infop1})
$$
Since~$V$ is $\cp$-quasi open, there exists an open set~$U$ contained
in~$\Om$ such that $\cp({U\bigtriangleup V},\Om)<\e$, and
by~(\rf{open}) there exists an open set~$W$ contained in~$\Om$ such
that ${U\bigtriangleup V}\subset W$ and $\cp(W,\Om)<\e$. As
$V\cup W=U\cup W$, the set $V\cup W$ is open. By subadditivity
(Theorem~\rf{subadd}) we have
$$
I\Le \ca(V\cup W,\Om)\Le \ca(V,\Om)+ \ca(W,\Om)\,.
$$
By~(\rf{upper}) and~(\rf{infop1}) we have
$$
I\Le \ca(V\cup W,\Om) \Le
\ca(E,\Om) + (1+k_1)\e + k_2(\Om) \e^{1/p}\,,
$$
where $k_1$ and $k_2(\Om)$ are the constants defined in~(\rf{const}).
Since ${\e>0}$ is arbitrary, we obtain $I\le \ca(E,\Om)$.
\endproof
\spazio
We conclude by proving that $\ca(\cdot,F)$ is countably subadditive.
\th{\thm{countsub}}{Let~$E$ and~$F$ be a bounded set in~$\rn$ and let
$(E_j)$ be a sequence of bounded sets in~$\rn$. If~$E$ is contained in
the union of the sequence $(E_j)$, then
$$
\ca(E,F) \Le \sum_{j=1}^\infty \ca(E_j,F) \,.
\leqno(\frm{countsub1})
$$
}
\negbigskip
\proof For every $k\in\n$ let $B_k=E\cap E_k$ and let
$G_k=B_1\cup\ldots\cup B_k$. By Theorems~\rf{monoE}
and~\rf{subadd} for every~$k$ we have
$$
\ca(G_k,F) \Le \sum_{j=1}^k \ca(B_j,F) \Le  \sum_{j=1}^\infty \ca(E_j,F)
\,.
$$
Since~$E$ is the union of the increasing sequence $(G_k)$, the continuity
along increasing sequences (Theorem~\rf{increasingE})
implies~(\rf{countsub1}).
\endproof
\parag{\chp{es}}{Capacity relative to a constant}
In this section we define the capacity $\ca(E,F,s)$ with respect to a
constant~$s$ by replacing the condition ${u=1}$ in $\partial E$ which
appears in~(\rf{dir0}) with the
condition ${u=s}$ in $\partial E$.
\defin{\thm{potentials}}{Let~$E$ and~$F$ be two $\cp$-compatible
bounded sets in~$\rn$, let~$s$ be a real number, and let $\psi$ be a
function in $\wp{\rn}$ such that ${\psi=1}$ \qe~$E$ and ${\psi=0}$
\qe~$F^c$. Every solution~$u$ of the
Dirichlet problem
$$
\dirpb{u-s\psi\in\wo{F\setminus E}}
{\A{u}=0}
{\wm{F\setminus E}},
\leqno(\frm{pots})
$$
is called a {\it $\ca$-potential\/} of~$E$ in~$F$ {\it relative to the
constant\/}~$s$.
The maximal and the minimal solutions of~(\rf{pots})
are called the {\it maximal\/} and
{\it minimal $\ca$-potentials\/} of~$E$ in~$F$ {\it relative to the
constant\/}~$s$.
}
\rem{\thm{potentials1}}{Clearly the previous definition does
not depend on the choice of~$\psi$. By the definition of the space
$\wo{F\setminus E}$ and by the properties of $\psi$ we have that
every $\ca$-potential~$u$ of~$E$ in~$F$ relative to the constant~$s$
satisfies $u=s$ \qe~$E$ and $u=0$ \qe~$F^c$.
}
\defin{\thm{caps}}{Let~$E$ and~$F$ be two $\cp$-compatible
bounded sets in~$\rn$ and let ${s\in\r}$. The
{\it capacity of~$E$ in~$F$ relative to the operator~$\A{}$ and to the
constant\/}~$s$ is defined as
$$
\ca(E,F,s)\Eq \langle \A{u},u\rangle \Eq
\iaa{F\setminus E}{Du}{Du}\,,
$$
where~$u$ is any $\ca$-potential of~$E$ in~$F$ relative to the
constant~$s$. By Corollary~\rf{b1}
this definition is independent of the choice of~$u$.
}
\rem{\thm{caps10}}{Let $s\neq0$ and let $a_s\colon\rn{\times}\rn\to
\rn$ be the function defined by
$$
a_s(x,\xi)=s\,\a{s\xi}\,.
$$
Then $a_s$ satisfies conditions (\rf{zero})--(\rf{bound}) with $c_1$,
$c_2$, and $b_2$ replaced by $|s|^{p}c_1$, $|s|^{p}c_2$, and
$|s|\,b_2$. Let $A_s$ be the operator defined by~(\rf{operator})
with~$a$
replaced by~$a_s$. Then~$u$ is a $\ca$-potential of~$E$ in~$F$ relative
to the constant~$s$ if and only if $u/s$ is a
$C_{A_s}$-potential of~$E$ in~$F$ according to Definition~\rf{potential},
and
$$
\ca(E,F,s)=C_{A_s}(E,F)\,.
$$
This shows that all properties of $\ca(E,F)$ proved in
Sections~\rf{capacity}--\rf{app} are still valid for $\ca(E,F,s)$, with some
obvious modifications, for every ${s\in\r}$.
}
\rem{\thm{caps10.5}}{From Lemma~\rf{signA} and Remark~\rf{caps10}
we obtain that, if~$u$ is a $\ca$-potential of~$E$ in~$F$ relative
to a constant ${s>0}$, then $\A{u}\ge0$ in $\wm{F}$ and $\A{u}\le0$ in
$\wm{E^c}$, whereas the inequalities are reversed when ${s<0}$. If
${s=0}$, by~(\rf{zero}) the function~$0$ is a
$\ca$-potential of~$E$ in~$F$ relative to the constant~$0$.
If~$u$ is another $\ca$-potential of~$E$ in~$F$ relative to the
constant~$0$, then $\A{u}=0$ in $\wm{\rn}$ by Theorem~\rf{b}.
}
\spazio
If $\A{u}=-\div\bigl(|Du|^{p-2}Du\bigr)$, then $\ca(E,F,s)=|s|^p\cp(E,F)$.
In the general case the relationship between $\ca(E,F,s)$ and $\cp(E,F)$
is given by the following proposition, which follows immediately from
Proposition~\rf{ca-cp} and Remark~\rf{caps10}.
\prop{\thm{caps11}}{Let~$E$ and~$F$ be two $\cp$-compatible
bounded sets in~$\rn$ and let ${s\in\r}$. Then
$$
\ldisplaylinesno{
\ca(E,F,s) \Ge |s|^p c_1\,\cp(E,F)\, -\, \|b_1\|_{\lone{F}} \,,
&(\frm{lowers})
\cr
\ca(E,F,s) \Le |s|^p k_1 \cp(E,F) \,+\, |s|\,k_2(F) \,\cp(E,F)^{1/p} \,,
&(\frm{uppers})
\cr}
$$
where $c_1$, $k_1$ and $k_2(F)$ are the constants which
appear in~(\rf{coerc}) and~(\rf{const}). If~$b_1$ and~$b_2$ belong to
$\linfty{F}$, then
$$
\ca(E,F,s) \Le \bigl(|s|^p k_1 + |s|\,k_3(F)\bigr) \,\cp(E,F) \,,
$$
where $k_3(F)$ is defined in~(\rf{const2}).
}
\spazio
The following lemma is an immediate consequence of the Comparison
Principle (Theorem~\rf{comparison}).
\lemma{\thm{caps14}}{Let~$E$ and~$F$ be two $\cp$-compatible
bounded sets in~$\rn$, let~$s_1$ and~$s_2$ be two real numbers, and
let~$u_1$ and~$u_2$ be the maximal (or minimal) $\ca$-potentials
of~$E$ in~$F$ relative to the constants~$s_1$ and~$s_2$ respectively. If
${s_1\le s_2}$, then ${u_1\le u_2}$ \qe~$\rn$.
}
\defin{\thm{caps14.5}}{Let~$E$ and~$F$ be two $\cp$-compatible
bounded sets in~$\rn$  and let ${s\in\r}$. We define
$$
\cas(E,F,s)=\cases{{1\over s}\ca(E,F,s)\,,
&if $s\neq0$,
\cr
\cr
0\,,
&if $s=0$.
\cr}
$$
}
\rem{\thm{caps14.7}}{By~(\rf{pots}) we have
$\cas(E,F,s)=\langle \A{u}, v\rangle$ for every $\ca$-potential of~$E$
in~$F$ relative to the constant~$s$ and for every function~$v$ in
$\wp{\rn}$ such that ${v=1}$ \qe~$E$ and ${v=0}$ \qe~$F^c$ (see
Remark~\rf{caps10.5} for the case ${s=0}$).
}
\spazio
We prove now that $\cas(E,F,s)$ depends continuously on~$s$.
\th{\thm{caps15}}{Let~$E$ and~$F$ be two $\cp$-compatible
bounded sets in~$\rn$. Then the function $s\mapsto\cas(E,F,s)$ is
continuous on~$\r$.
}
\proof Let $\r_+=\{{s\in\r:s\ge0}\}$ and $\r_-=\{{s\in\r:s\le0}\}$. We
prove only that ${s\mapsto\cas(E,F,s)}$ is
continuous on~$\r_+$, the proof for~$\r_-$ being analogous. We begin by
proving the right continuity on~$\r_+$.
Let us fix $s\ge0$ and let $(s_j)$ be a decreasing sequence
in~$\r$ converging to~$s$. Let~$u$ and~$u_j$ be the maximal
$\ca$-potentials of~$E$ in~$F$ relative to the constants~$s$ and~$s_j$
respectively. As $\cp(E,F)\finite$, by~(\rf{uppers}) and~(\rf{coerc1}) the
sequence $(u_j)$ is bounded in $\wp{\rn}$, and by~(\rf{bound1}) the
sequence $(\A{u_j})$ is bounded in $\wm{\rn}$. Passing, if necessary, to
a subsequence, we may assume  that $(u_j)$ converges weakly in
$\wp{\rn}$ to some function $w\in\wo{F}$ and  that $(\A{u_j})$
converges weakly in $\wm{\rn}$ to some element~$f$ of $\wm{\rn}$.
By Lemma~\rf{caps14} we have $u\le u_j\le u_i$ \qe~$\rn$ for
every ${j\ge i}$, hence $u\le w\le u_i$ \qe~$\rn$. Since ${u=s}$ and
${u_i=s_i}$ \qe~$E$ (Remark~\rf{potentials1}), as ${i\to\infty}$ we
obtain that ${w=s}$ \qe~$E$.

We want to prove that $f=\A{w}$ in $\wm{F}$ and that ${w=u}$
\qe~$\rn$. From the monotonicity condition~(\rf{mono1}) for every~$j$
we
obtain
$$
\langle \A{v}, v-u_j \rangle \Ge \langle \A{u_j}, v-u_j \rangle
\qquad \forall v\in\wp{\rn}\,.
\leqno(\frm{100})
$$
If ${j\ge i}$, by Lemma~\rf{caps14} we have
$u_j\le u_i$ \qe~$\rn$. Since $\A{u_j}\ge0$ in $\wm{F}$
(Remark~\rf{caps10.5}), we have
$$
\langle \A{u_j}, v-u_j \rangle
\Ge \langle \A{u_j}, v-u_i \rangle \qquad\forall v\in\wo{F}\,,
\leqno(\frm{101})
$$
which, together with~(\rf{100}), gives
$$
\langle \A{v}, v-u_j \rangle \Ge \langle \A{u_j}, v-u_i \rangle
\qquad \forall v\in\wo{F}
$$
whenever ${j\ge i}$. Passing to the limit as ${j\to\infty}$ we obtain
$$
\langle \A{v}, v-w \rangle \Ge \langle f , v-u_i \rangle
\qquad \forall v\in\wo{F}\,,
$$
and as ${i\to\infty}$ we get
$$
\langle \A{v}, v-w \rangle \Ge \langle f , v-w\rangle
\qquad \forall v\in\wo{F}\,.
\leqno(\frm{102})
$$
Putting $v=w+\e z$, with $z\in\wo{F}$ and ${\e>0}$, and dividing
by~$\e$ we obtain
$$
\langle A(w+\e z), z \rangle \Ge \langle f , z\rangle
\qquad\forall z\in \wo{F}\,.
$$
Passing to the limit as $\e\to0$ we get
$$
\langle \A{w}, z\rangle \Ge \langle f , z\rangle
\qquad\forall z \in \wo{F}\,,
$$
hence $\A{w}=f$ in $\wm{F}$.

By the definition of~$u_j$ we have
$\langle \A{u_j}, v \rangle =0$  for every
$v\in\wo{{F\setminus E}}$ and for every~$j$.
As ${j\to\infty}$ we obtain
$$
\langle \A{w}, v \rangle \Eq \langle f, v \rangle \Eq
\lim_{j\to\infty} \langle \A{u_j}, v \rangle \Eq0
\qquad \forall v\in\wo{F\setminus E}\,,
$$
hence~$w$ is a $\ca$-potential of~$E$ in~$F$ relative to the
constant~$s$. Since ${u\le w}$ \qe~$\rn$, by the maximality of~$u$ we
obtain ${u=w}$ \qe~$\rn$.

Let~$\psi$ be a function in $\wp{\rn}$ such that ${\psi=1}$ \qe~$E$ and
${\psi=0}$ \qe~$F^c$.
Since $(\A{u_j})$ converges to~$f$ weakly in $\wm{\rn}$ and $\A{u}=f$
in $\wm{F}$, by Remark~\rf{caps14.7} we have
$$
\cas(E,F,s)\Eq \langle \A{u},\psi\rangle \Eq \langle f,\psi\rangle
\Eq \lim_{j\to\infty} \langle \A{u_j},\psi\rangle
\Eq \lim_{j\to\infty} \cas(E,F,s_j) \,.
$$
This proves that the function ${s\mapsto\cas(E,F)}$ is right continuous
on~$\r_+$. For the proof of the left continuity on~$\r_+$,
we fix ${s>0}$ and an increasing sequence $(s_j)$ converging to~$s$. We
may assume that ${s_j>0}$ for every~$j$. Then we use the same
arguments as in the first part of the proof, with the only difference that
now we use the minimal $\ca$-potentials instead of the maximal
$\ca$-potentials. As $(s_j)$ is increasing, the sequence $(u_j)$ is
increasing, and, consequently, $u_i$ must be replaced by~$w$
in~(\rf{101}) and we obtain directly~(\rf{102}). The final part of the
proof remains unchanged.
\endproof
\spazio
We prove now that $\cas(E,F,s)$ is increasing with
respect to~$s$.
\th{\thm{monos}}{Let~$E$ and~$F$ be two  $\cp$-compatible
bounded sets in~$\rn$ and let $s_1$ and $s_2$ be two real numbers
with ${s_1<s_2}$. Then $\cas(E,F,s_1) \le \cas(E,F,s_2)$.
}
\proof Let~$u_1$ and~$u_2$ be two $\ca$-potentials of~$E$ in~$F$
relative to the constants~$s_1$ and~$s_2$ respectively, let~$t$ be a real
number such that $0<t<s_2-s_1$, and let~$v$ be the function of
$\wp{\rn}$ defined by
$v={1\over t}\bigl((u_2-u_1)\land t\bigr)$. By Remark~\rf{potentials1}
we have $v=1$ \qe~$E$ and $v=0$ \qe~$F^c$.
Therefore, Remark~\rf{caps14.7} implies that
$$
\displaylines{
\cas(E,F,s_2) - \cas(E,F,s_1) \Eq
\langle \A{u_2}-\A{u_1}, v\rangle \Eq
\cr
\Eq {1\over t}  \iaaa {\{u_2-u_1<t\}} {Du_2} {Du_1} {Du_2-Du_1} \,,
\cr}
$$
and the conclusion follows from the monotonicity
condition~(\rf{mono}).
\endproof
\intro{Acknowledgments}
This work is part of the Project EURHomogenization,
Contract SC1-CT91-0732 of the Program SCIENCE of the Commission of
the European  Communities, and of the Research Project ``Irregular
Variational Problems'' of the Italian National Research Council.
\intro{References}
\def\interrefspace{\smallskip}  
{\ninepoint\frenchspacing
\item{[\bib{CAS-GAR}]}CASADO DIAZ  J., GARRONI  A.: A compactness
theorem for Dirichlet problems in varying domains with monotone
operators. In preparation.
\interrefspace
\item{[\bib{CHO}]}CHOQUET  G.: Forme abstraite du th\'eor\`eme de
capacitabilit\'e. {\it Ann. Inst. Fourier (Grenoble)\/} {\bf 9} (1959),
83-89.
\interrefspace
\item{[\bib{DM-83}]}DAL MASO  G.: On the integral representation of
certain local functionals. {\it Ricerche Mat.\/} {\bf 32} (1983), 85-113.
\interrefspace
\item{[\bib{DM-DEF}]}DAL MASO  G., DEFRANCESCHI  A.: Limits of
nonlinear Dirichlet problems in varying domains. {\it Manuscripta
Math.\/} {\bf 61} (1988), 251-278.
\interrefspace
\item{[\bib{DM-GAR-cap}]}DAL MASO  G., GARRONI  A.: Capacity theory
for non-symmetric elliptic operators. Preprint SISSA, Trieste, 1993.
\interrefspace
\item{[\bib{DM-GAR}]}DAL MASO  G., GARRONI  A.: New results on the
asymptotic behaviour of Dirichlet problems in perforated domains. {\it
Math. Mod. Meth. Appl. Sci.\/} {\bf 3} (1994), 373-407.
\interrefspace
\item{[\bib{DM-GAR-94}]}DAL MASO  G.,  GARRONI  A.: The capacity
method for asymptotic Dirichlet problems. Preprint SISSA, Trieste, 1994.
\interrefspace
\item{[\bib{DM-MUR-lum}]}DAL MASO  G., MURAT  F.: Dirichlet
problems in perforated domains for homogeneous monotone operators
on $H^1_0$. {\it Calculus of Variations, Homogenization and Continuum
Mechanics (CIRM-Luminy, Marseille, 1993)\/}, 177-202, {\it World
Scientific, Singapore\/}, 1994.
\interrefspace
\item{[\bib{DM-MUR}]}DAL MASO  G., MURAT  F.: Asymptotic behaviour
and correctors for Dirichlet problems in perforated domains with
homogeneous monotone operators. In preparation.
\interrefspace
\item{[\bib{EVA-GAR}]}EVANS L.C., GARIEPY R.F.: Measure Theory and
Fine Properties of Functions. CRC Press, Boca Raton, 1992.
\interrefspace
\item{[\bib{FED-ZIE}]}FEDERER  H., ZIEMER  W.P.: The Lebesgue set of a
function whose distribution derivatives are $p$-th power summable.
{\it Indiana Univ. Math. J.\/} {\bf 22} (1972), 139-158.
\interrefspace
\item{[\bib{HEI-KIL-MAR}]}HEINONEN  J., KILPEL\"AINEN T., MARTIO
O.:
Nonlinear Potential Theory of Degenerate Elliptic Equations. Clarendon
Press, Oxford, 1993.
\interrefspace
\item{[\bib{KIN-STA}]}KINDERLEHRER  D., STAMPACCHIA  G.: An
Introduction to Variational Inequalities and Their Applications.
Academic Press, New York, 1980.
\interrefspace
\item{[\bib{MAZ}]} MAZ'YA  V.G.: Sobolev Spaces. Springer-Verlag,
Berlin, 1985.
\interrefspace
\item{[\bib{SKR-86}]}SKRYPNIK  I.V.: Nonlinear Elliptic Boundary Value
Problems. Teubner-Verlag, Leipzig, 1986.
\interrefspace
\item{[\bib{SKR-90}]}SKRYPNIK  I.V.: Methods of Investigation of
Nonlinear Elliptic Boundary Value Problems (in Russian). Nauka,
Moscow, 1990.
\interrefspace
\item{[\bib{SKR-91}]}SKRYPNIK  I.V.: Averaging nonlinear Dirichlet
problems in domains with channels. {\it Soviet Math. Dokl.\/} {\bf 42}
(1991), 853-857.
\interrefspace
\item{[\bib{SKR-93}]}SKRYPNIK  I.V.: Asymptotic behaviour of solutions
of nonlinear elliptic problems in perforated domains. {\it Mat. Sb.
(N.S.)\/} {\bf 184} (1993), 67-90.
\interrefspace
\item{[\bib{SKR-94}]}SKRYPNIK  I.V.: Homogenization of nonlinear
Dirichlet problems in perforated domains of general structure. Preprint
SISSA, Trieste, 1994.
\interrefspace
\item{[\bib{STA}]}STAMPACCHIA  G.: Le probl\`eme de Dirichlet pour
les \'equations elliptiques du second ordre \`a coefficients discontinus.
{\it Ann. Inst. Fourier (Grenoble)\/} {\bf 15} (1965), 189-258.
\interrefspace
\item{[\bib{ZIE}]}ZIEMER  W.P.: Weakly Differentiable Functions.
Springer-Verlag, Berlin, 1989.
\par
}
\bye